\documentclass[a4paper,graphics,floatfix,nofootinbib,tightenlines,nobibnotes,arp,11pt]{revtex4}
\usepackage{lipsum}
\usepackage[leqno]{amsmath}
\usepackage[english]{babel}
\usepackage{fancyhdr}
\usepackage[usenames, dvipsnames]{color}
 
\definecolor{mypink1}{rgb}{0.858, 0.188, 0.478}
\definecolor{mypink2}{RGB}{219, 48, 122}
\definecolor{mypink3}{cmyk}{0, 0.7808, 0.4429, 0.1412}
\definecolor{mygray}{gray}{0.6}
\definecolor{Myorange}{cmyk}{0,0.42,1,0}
\definecolor{Myblue}{rgb}{0.1,0.1,1}
\definecolor{Mygrey}{gray}{0.75}
\definecolor{Mygreen}{named}{SpringGreen}
\usepackage{hyperref}
\usepackage{calc}
\usepackage{newlfont}
\usepackage[english]{babel}
\usepackage{amssymb,amsmath,amsfonts}
\usepackage{mathrsfs}
\usepackage{dsfont}
\usepackage[all]{xy}
\usepackage{enumerate}
\usepackage{textcomp}
\usepackage{pifont}
\usepackage{amsthm}
\usepackage{bm}
 \usepackage{lineno}
\DeclareMathOperator{\Char}{Char}
\DeclareMathOperator{\WF}{WF}

\newtheorem{theorem}{Theorem}
\newtheorem{proposition}{Proposition}
\newtheorem{lemma}{Lemma}
\newtheorem{corollary}{Corollary}
\newtheorem{definition}{Definition}

\newcounter{obsctr}

\newtheorem{remark}{Remark}
\makeatletter
\newcommand{\subjclass}[2][2010]{%
  \let\@oldtitle\@title%
  \gdef\@title{\@oldtitle\footnotetext{#1 \emph{Mathematics subject classification.} #2}}%
}
\renewcommand{\keywords}[1]{%
  \let\@@oldtitle\@title%
  \gdef\@title{\@@oldtitle\footnotetext{\emph{Key words and phrases.} #1.}}%
}
\makeatother
\newcommand\shorttitle{Gevrey regularity for Sums of Squares}
\newcommand\authors{$\mathscr{G}\!.$ \textit{Chinni}}

\begin{document}
%
\title{On the Gevrey regularity for Sums of Squares of vector fields,\\ study of some models.}
\vspace*{2em}
\author{G. Chinni}
\email[\ding{41}  ]{gregorio.chinni@gmail.com}
\altaffiliation{{\it Fakult\"at f\"ur Mathematik, Oskar--Morgenstern--Platz 1, 1090 Vienna, Austria}}
\keywords{Sums of squares, Microlocal Gevrey regularity, Hypoellipticity\vspace{0.2em}}
\subjclass[2010]{35H10, 35H20, 35B65.\vspace{0.3em}}
\maketitle
\tableofcontents
\vspace{1.5cm}
\abstractname{.
The local and micro-local Gevrey hypoellipticity of a class of ``sums of squares" 
with real analytic coefficients is studied in detail.
Some partial regularity result is also given.}
\newpage
\section{Introduction}
The purpose of this paper is to discuss the Gevrey hypoellipticity
properties of three model operators  that are sums of squares of
vector fields in four dimensions. The operators have analytic coefficients
and verify the H\"ormander condition: 
the Lie algebra generated by the vector fields as well as by their commutators
has, in every point, dimension equal to the dimension of the ambient space.
Hence in view of the celebrated H\"ormander theorem, \cite{H67},
the operators are $C^{\infty}$-hypoelliptic.\\
Let $P(x;D)= \sum_{1}^{k}X_{j}^{2}(x,D)$, $X_{j}(x,D)$ vector fields with real 
analytic coefficients on $ \Omega $ open subset $\mathbb{R}^{n}$.
We say that $ P $ is   $C^{\infty}\left( G^{r}\right)$-hypoelliptic,
$ r \geq 1$, in $\Omega$ if for every $U$ open subset of $\Omega$
and every $ u \in \mathscr{D}'\left( U \right)$, $Pu \in C^{\infty}(U)$ 
$ \left( G^{r}(U)\right)$ implies $ u \in C^{\infty}(U)$ 
$\left( G^{r}(U)\right)$.
When $ r=1 $ we say that $P$ is analytic hypoelliptic. We recall that
$G^{r}(U)$ denotes the $r-$Gevrey class of function on $U$:
a $C^{\infty}-$function, $f$, on $U$ belongs to $G^{r}(U)$,
$1\leq r\leq \infty$, if for every $K$ compact subset of $U$
there is a constant $C_{K}$ such that 
$\displaystyle|\partial^{\alpha}f(x)|\leq C_{K}^{|\alpha|+1}(\alpha!)^{r}$ 
for every $\alpha \in \mathbb{N}^{n}$ and $x \in K$.\\
Derridj showed in \citep{D-71} that  for $P$ as above,
the H\"ormander condition is necessary for the analytic hypoellipticity but it is not sufficient.
An example of operator sum of squares of real analytic vector fields
satisfying the H\"ormander condition but not analytic hypoelliptic
was given by Baouendi and Goulaouic in \citep{BG-72}.
At the present, there aren't general analytic hypoellipticity results.
Some results, in this direction, were obtained
by Treves, \cite{Tr-78}, Tartakoff, \cite{Tart-80}, and Albano and Bove, 
\citep{AB-13}.
For completeness, we recall that, with regard to the Gevrey regularity, if no 
additional assumption is made on the operator $P$, the (local) optimal 
characterization was obtained by Derridj and Zuily, \cite{DZ}.
In 1999 Treves formulated a conjecture which related the analytic hypoellipticity with 
geometrical properties of the characteristic variety of $P$, see \cite{Tr-99} and \cite{Tr-2006}.\\
In recent papers Albano, Bove and Mughetti, \cite{ABM-TrC1-2016}, and
Bove and Mughetti, \cite{Bove_Mughetti-TrC2-2016},
showed that the sufficient part of the Treves' conjecture does not hold
neither locally nor microlocally.
More precisely in \cite{ABM-TrC1-2016} and \cite{Bove_Mughetti-TrC2-2016}
the authors produced and studied the first models which are not consistent with
the Treves conjecture, \cite{Tr-2006}. However, contrary to the cases of 
\cite{ABM-TrC1-2016} and \cite{Bove_Mughetti-TrC2-2016},
the operators studied here have no exceptional strata because
the symbols do not depend on the tangent variables of the
``inner most" stratum.\\
Our results can be stated as follows:
\begin{theorem}\label{T-1}
Let $P_{1}(x,D)$ the sum of squares given by
\vspace*{-0.4em}
\begin{equation}\label{Op_1} 
D_{1}^{2} + x_{1}^{2(p-1)} D_{2}^{2} + x_{1}^{2(q-1)} D_{3}^{2}
+ x_{1}^{2(r-1)} x_{2}^{2k} D_{4}^{2} + x_{1}^{2(r+ \ell -1 )} D_{4}^{2},
\end{equation}
where $ p, q, r, k $ and $ \ell $ are positive integers such that $ p < q < r $
and $pk < \ell $ and $P_{2}(x,D)$ the sum of squares given by
\vspace*{-0.4em}
\begin{equation}\label{Op_2}
D_{1}^{2} + x_{1}^{2(p-1)} D_{2}^{2} + x_{1}^{2(q-1)} D_{3}^{2}
+ x_{1}^{2(r-1)} x_{3}^{2k} D_{4}^{2} + x_{1}^{2(r+ \ell -1 )} D_{4}^{2},
\end{equation}
where $ p, q, r, k $ and $ \ell $ are positive integers such that $ p < q < r $
and $qk < \ell $.
We have:
\begin{itemize}
\vspace*{-0.4em}
\item[\textbf{i)}]  $P_{1}(x,D)$ is $G^{s}$-hypoelliptic with 
             $s=\sup\Big\lbrace\frac{r+kp}{q}, \frac{r}{p}\Big\rbrace$.
\vspace*{-0.4em}             
\item[\textbf{ii)}] $P_{2}(x,D)$ is $G^{s}\!$-hypoelliptic with 
              $s\!= \frac{r+kq}{p}$. 
\end{itemize}
\end{theorem}
\noindent
The strategy used to obtain the above results shows, without particular 
technical trouble, that:
%
\begin{remark}\label{R-T_1}
If $pk \geq \ell $ then $P_{1}$ is $G^{s}$-hypoelliptic with
$s = \sup\Big\lbrace\frac{r+\ell}{q}, \frac{r}{p}\Big\rbrace$
and $P_{2}$ is $G^{s}$-hypoelliptic with $s = \frac{r+\ell}{p}$.
\end{remark}
%
We recall that by the result of Derridj and Zuily, \cite{DZ}, 
$P_{1}$ is $(r+kp)$-Gevrey hypoelliptic and $P_{2}$ is $(r+kq)$-Gevrey 
hypoelliptic when $kq < \ell$ and they are both $(r+\ell)$-Gevrey hypoelliptic
when $kp \geq \ell$. 
%
%
\begin{theorem}\label{T-2}
Let the operator $P_{3}(x,D)$ be given by
\begin{equation}\label{Op_3}
D_{1}^{2} + x_{1}^{2(p-1)} D_{2}^{2} + x_{1}^{2(q-1)} D_{3}^{2}
+ x_{1}^{2(r-1)} x_{2}^{2k} D_{4}^{2} + x_{1}^{2(f -1 )}x_{3}^{2\ell}
D_{4}^{2} + x_{1}^{2(f+e-1)}D_{4}^{2}
\end{equation}
in $ \Omega $, open neighborhood of the origin in $ \mathbb{R}^{4} $,
where $ p, q, r, k, f, \ell $ and $ e $ are positive integers such that 
$ p < q < r < f$ and $e >\sup\lbrace pk,q\ell\rbrace$.
We have:
$P_{3}$ is $G^{s}$-hypoelliptic, with
$s = \displaystyle\sup\Big\lbrace\frac{r+kp}{q}, \frac{r}{p}\Big\rbrace$
if $f > r + kp$ and
$s = \displaystyle\sup\Big\lbrace\frac{f}{q}, \frac{r}{p}\Big\rbrace$
if $f \leq r + kp$.
\end{theorem}
%
We point out that, in accordance with the results in \cite{DZ},
the operator $P_{3}$
is $(r+kp)$-Gevrey hypoelliptic if $f > r+kp $ and 
$\inf\lbrace r+kp, f+q\ell\rbrace$-Gevrey hypoelliptic if $f <r+kp$.\\
The strategy used to proof the Theorem \ref{T-2} shows, without particular 
technical trouble, that:
\begin{remark}
If $e <\sup\lbrace pk,q\ell\rbrace$ we can distinguish two cases:
$f + e < r+ kp$, the operator $P_{3}$ is a generalization of the 
Ole\u{\i}nik-Radkevi\v c  operator and it is $G^{\frac{f+e}{p}}$-hypoelliptic 
and $ kp < e < lq$, i. e. $r+kp < f+e$, $P_{3}$ 
is a generalization of the operator $P_{2}$, it is $G^{s}$-hypoelliptic
with $s=\sup\lbrace \frac{r+kp}{q}, \frac{r}{p}\rbrace$.
\end{remark}
\begin{remark}
If $p<q<f<r$ and $e >\sup\lbrace r+pk,f+q\ell\rbrace$
we can distinguish two cases:
$f +q\ell > r $, $P_{3}$ is $G^{\frac{r}{p}}$-hypoelliptic,
and $f +q\ell < r $, $P_{3}$ is $G^{\frac{f+q\ell}{p}}$-hypoelliptic. 
\end{remark}
\noindent
Even if, at the present, the proof of the optimality of the operators
$P_{1}$, $P_{2}$ and $P_{3}$ is an open pro\-blem, we think that the Gevrey 
regularities obtained are optimal.

\vspace{0.3cm}
\begin{remark}
\vspace*{-0.7em}
The results stated above can be extended to the operators
\begin{equation}\nonumber
P(x,D)=D_{1}^{2} +\sum_{j=2}^{n-1} x_{1}^{2(r_{j}-1)} D_{j}^{2} 
+ \left( x_{1}^{2(r_{n}-1)} x_{i}^{2k} + x_{1}^{2(r_{n}+ \ell -1 )} 
\right)D_{n}^{2}, \qquad 2\leq i\leq n-1,
\end{equation}
defined in $\Omega$, open neighborhood of the origin in $\mathbb{R}^{n}$,
where $ r_{j}$, $j= 1,\dots, n$, $k$ and $\ell$ are positive integers such that
$r_{1}< r_{2} < \dots < r_{n}$.
We have: if $kr_{i} <\ell$, $P(x,D)$ is $G^{s}$-hypoelliptic with 
$s= \displaystyle\sup\Big\lbrace\frac{r_{n}+kr_{2}}{r_{3}}, \frac{r_{n}}
{r_{2}}\Big\rbrace$ if $i =2$ and $s= \displaystyle\frac{r_{n}+kr_{i}}{r_{2}}$ 
if $i \neq 2$; if $kr_{i} \geq\ell$, $P(x,D)$ is $G^{s}$-hypoelliptic with 
$s= \displaystyle\sup\Big\lbrace\frac{r_{n}+\ell}{r_{3}}, \frac{r_{n}}{r_{2}}
\Big\rbrace$ if $i =2$ and $s= \displaystyle\frac{r_{n}+\ell}{r_{2}}$
if $i \neq 2$.
This situation does not present additional difficulties compared
to that we are going to handle.
\end{remark}

In the last section we analyze the partial regularity, that is, following the ideas in 
\cite{BT} , we study the non-isotropic Gevrey regularity of $P_1$ and $P_2$. 
In addiction to give a more precise characterization of the local regularity of the 
operators in a neighborhood of the origin, the purpose is to make in evidence 
further differences with the Ole\u{\i}nik-Radkevi\v c  operator: as shown the 
operator $P_{1}$, both in the case $pk< \ell$ that in the case $pk\geq \ell$, 
does not have directions with analytic growth, thing which, on the contrary,
occurs both in $P_{2}$ and  in the Ole\u{\i}nik-Radkevi\v c  operator.
\bigskip
\paragraph*{Acknowledgements:}The author was partially supported by a postdoctoral fellowship from FAPESP,
Grant 2013/08238-6.
The author is also happy to thank the ``Instituto de Matem\'atica e Estat\'{\i}stica'' of the University  of S\~ao Paulo  for its generous and kind hospitality and in particular the research group on Partial Differential Equations and Complex Analysis for the excellent working conditions.
%
\section{Proof of Theorem \ref{T-1}}\label{P-T_1} 
\subsection{Gevrey Regularity for $P_{1}(x,D)$.}
\noindent
The characteristic variety of $P_{1}$ is:
\begin{equation}\nonumber
Char(P_{1}) =\lbrace (x;\xi)\in T^{*}\mathbb{R}^{4}\setminus \lbrace 0 \rbrace\,
:\, x_{1}=0,\, \xi_{1}=0\rbrace.
\end{equation}
Following the ideas in \cite{Tr-2006} and \cite{BoTr} it can be seen as the disjoint union of analytic submanifolds, strata: 
\begin{equation}\nonumber
Char(P_{1}) = \Sigma_{p,\pm} \cup \Sigma_{q,\pm} \cup \Sigma_{r,\pm}
\cup \Sigma_{r+kp,\pm},
\end{equation}
where
\begin{equation}\nonumber
\begin{array}{rcl}
\Sigma_{p,\pm}&=&
\lbrace  (x;\xi) \in T^{*}\mathbb{R}^{4}\setminus \lbrace 0
\rbrace :\, x_{1}=0,\, \xi_{1}=0, \,\xi_{2}\gtrless 0 \rbrace,
\\[5pt]
\Sigma_{q,\pm}  &=& 
\lbrace  (x;\xi)\in T^{*}\mathbb{R}^{4}\setminus \lbrace 0
\rbrace :\, x_{1}=0,\, \xi_{1}=0,\, \xi_{2} =0, \,\xi_{3}\gtrless 0 \rbrace,
\\[5pt]
\Sigma_{r,\pm} &=& 
\lbrace  (x;\xi)\in T^{*}\mathbb{R}^{4}\setminus \lbrace 0
\rbrace :\, x_{1}=0,\, \xi_{1}=0,\, \xi_{2} =0,\, \xi_{3}=0,
\,\xi_{4}x_{2}\neq 0 
\rbrace
\\[5pt]
\Sigma_{r+pk,\pm} &=&  
\lbrace  (x;\xi)\in T^{*}\mathbb{R}^{4}\setminus \lbrace 0 
\rbrace :\,
x_{1}=0,\, \xi_{1}=0,\, \xi_{2} =0,\, \xi_{3}=0,\, x_{2}= 0,\, \xi_{4}\gtrless 0 \rbrace.
\end{array}
\end{equation}
\noindent
By the results of Derridj and Zuily, \cite{DZ}, and Roth\-schild and Stein, 
\cite{RS}, the operator $ P_{1}$ has the following sub-elliptic 
estimate with loss of $2(1-1/(r+kp))$ derivatives:
%
\begin{equation}\label{SB-E-1/2}
\|u\|^{2}_{\frac{1}{r+kp}} + \sum_{j=1}^{5} \|X_{j} u\|^{2} 
\leq 
C \big( |\langle P_{1}u, u \rangle | + \| u \|^{2} \big). 
\end{equation}
Here $ X_{1} = D_{1} $, $ X_{2} = x_{1}^{p-1}D_{2}$,
$ X_{3} = x_{1}^{q-1} D_{3}$, $ X_{4} = x_{1}^{r-1} x_{2}^{k} D_{4} $,
$ X_{5} = x_{1}^{r+\ell-1} D_{4}$,
$ \| \cdot \|_{s} $ denotes the $ H^{s} $ Sobolev norm and 
$ \| \cdot \| = \| \cdot \|_{0} $ denotes the $ L^{2} $ norm on the fixed
open set $ \Omega $.\\   
To study the regularity of the solutions we estimate 
the high order derivatives of the solutions in $ L^{2} $ norm.
As a matter of fact we estimate a suitable localization of a high derivative 
using the above estimate.
For $ x_{1} \neq 0 $ the operator $ P_{1} $ is elliptic and we shall not examine
this region, elliptic operators are Gevrey hypoelliptic in any class $ G^{s} $ for
$ s \geq 1 $.\\
Let $ \varphi_{N} (x) $ be a cutoff function of Ehrenpreis-H\"or\-man\-der
type: $\varphi_{N} \in C^{\infty}_{0}\left( \Omega \right)$ non negative
such that $ \varphi_{N} \equiv 1 $ on $ \Omega_{0} $, 
$ \Omega_{0} $ neighborhood of the origin compactly 
contained in $ \Omega $, and exists a constant $ C $ such that for every 
$ | \alpha | \leq 2 ( r + pk ) N $, $ \alpha \in \mathbb{N}^{4} $, we have 
$ | D^{\alpha} \varphi_{N} | \leq C^{ |\alpha| +1} N^{| \alpha | }$.

\noindent
We may assume that $ \varphi_{N} $ is independent of  the $ x_{1} $-variable
since every $ x_{1} $-derivative landing on $ \varphi_{N} $ would leave a
cut off function supported where $ x_{1} $ is bounded away from zero,
where the operator is elliptic. Moreover we may assume that $ \varphi_{N} $
is independent of the $ x_{2} $-variable since every $ x_{2}$-derivative
landing on $ \varphi_{N} $ would leave a cut off function supported where
$ x_{2} $ is bounded away from zero, in this region the operator 
satisfies the H\"or\-man\-der-Lie algebra condition at the step $ r $.
The operator $ P_{1} $, in this region, has the following estimate with loss of
$ 2 ( 1 - 1/r ) $ derivatives:
\begin{equation}\nonumber
\|u\|^{2}_{\frac{1}{r}} + \sum_{j=1}^{5} \|X_{j} u\|^{2}
\leq 
C \big( |\langle Pu, u \rangle | + \| u \|^{2} \big),
\end{equation}
where $ u \in C^{\infty}_{0} ( K )$ with 
$ K \cap \lbrace x_{2}=0  \rbrace = \emptyset $.
In this region the operator is a generalization of the Ole\u{\i}nik-Radkevi\v c
operator then $ P_{1}$ is $ G^{r/p}$-hypoelliptic and not better,
for more details see \cite{C} and \cite{BT}.
Then, we can conclude that 
if $v$ solves the equation $P_{1}v = f$ and $f$ is analytic then the points
$\displaystyle\rho_{_{1}}\in \Sigma_{r}$ does not belong
to $WF_{\frac{r}{p}}(v)$.\\
Now, we are interested to the microlocal regularity 
in  $\displaystyle\rho_{_{0}}\in \Sigma_{r+kp}$. To obtain
this it is sufficient to study the microlocal regularity of $P_{1}$ in 
$(0; e_{4})$. Indeed
the microlocal regularity in  a generic point $\displaystyle\rho_{_{0}}$ can be 
obtained following the same strategy below with the only difference that the 
cut-off function $\varphi_{N}(x)$ will be identically equal to 1 in $\Omega_{0}$ 
neighborhood of $\pi_{x}(\rho_{_{0}})=(0,0,x_{3}^{0},x_{4}^{0})$,
where $\pi_{x}$ is the projection in the space variables.
Thus since we are interested to the microlocal regularity of $P_{1}$
in $(0 ; e_{4})$ we take $\varphi_{N} (x) = \varphi_{N}(x_{3},x_{4})$.\\
%
We replace $ u $ by $ \varphi_{N} (x) D_{4}^{N} u $ in ( \ref{SB-E-1/2}).
We have
\begin{equation}\label{SB-E-1.1_D4}
\| \varphi_{N} D^{N}_{4} u\|^{2}_{\frac{1}{r+kp}}\! 
+\! \sum_{j=1}^{5} \|X_{j} 
\varphi_{N} D^{N}_{4} u \|^{2}
\leq \!
C 
\left( |\langle P_{1} \varphi_{N} D^{N}_{4} u, \varphi_{N} D^{N}_{4} u\rangle | 
+ \| \varphi_{N} D^{N}_{4} u \|^{2}  \right). 
\end{equation}
The scalar product in the right hand side leads to
\begin{equation}\label{BR-1.1}
\begin{split} 
& \langle  \varphi_{N} D_{4}^{N} P_{1}  u, \varphi_{N} D_{4}^{N} u \rangle
+ \sum_{j=1}^{5} \langle \lbrack X_{j}^{2}, \varphi_{N} D_{4}^{N} \rbrack u, 
\varphi_{N} D_{4}^{N}u \rangle
\\
&\qquad\quad
= 2 \sum_{j=1}^{5}
 \langle  \lbrack X_{j}, \varphi_{N} D_{4}^{N} \rbrack u,
 X_{j}\varphi_{N} D_{4}^{N}u \rangle
+\sum_{j=1}^{5} 
 \langle \lbrack \lbrack X_{j}, \varphi_{N} D_{4}^{N} \rbrack , 
 X_{j}\rbrack u, \varphi_{N} D_{4}^{N} u \rangle
\\
&\qquad\qquad\qquad \qquad\qquad\qquad\qquad\qquad\qquad\qquad
\qquad
\qquad
+  \langle  \varphi_{N} D_{4}^{N} P_{1}  u, \varphi_{N} D_{4}^{N} u \rangle.
\end{split}
\end{equation}
The last term is trivial to estimate since $ P_{1} u $ is analytic; we may assume 
without loss of generality, that is zero. Since $ \varphi_{N} $ 
depends only by $ x_{3} $ and $ x_{4} $
we must analyze the commutators with $ X_{3}$, $ X_{4}$ and $ X_{5}$.
Before to give the general form of the terms which appear inside of the iterating 
process we begin to analyze some particular situations.\\
Case $ X_{4} $. We have
\begin{equation}\label{es-X_4}
\begin{split}
&
| \langle  \lbrack X_{4}, \varphi_{N} D_{4}^{N} \rbrack u,
X_{4}\varphi_{N} D_{4}^{N}u \rangle|
+ 
| \langle \lbrack \lbrack X_{4}, \varphi_{N} D_{4}^{N} \rbrack , X_{4}\rbrack u, 
\varphi_{N} D_{4}^{N} u \rangle| 
\\
&\quad
= 2 | \langle x_{1}^{r-1}x_{2}^{k} \varphi^{(1)}_{N} D_{4}^{N} u,
X_{4}\varphi_{N} D_{4}^{N}u \rangle |
+ | \langle  x_{1}^{2(r-1)} x_{2}^{2k} \varphi^{(2)}_{N} D_{4}^{N} u,
\varphi_{N} D_{4}^{N} u \rangle |.
\end{split}
\end{equation}
The first term, we have
\begin{equation}\label{es-X_4-1}
\begin{split}
| \langle x_{1}^{r-1}x_{2}^{k} \varphi^{(1)}_{N} D_{4}^{N} u,
X_{4}\varphi_{N} D_{4}^{N}u \rangle |
&\leq  \sum_{j=1}^{ N}  C_{j} \| X_{4} \varphi^{(j)}_{N} D_{4}^{N - j} u \|^{2}
\\
&\qquad
+  \sum_{j=1}^{ N+1}\!\frac{1}{C_{j}} \| X_{4}\varphi_{N} D_{4}^{N}u \|^{2}
+  C_{N+1}\| \varphi^{(N+1)}_{N}  u \|^{2} ,
\end{split}
\end{equation}
the constants $ C_{j}$ are arbitrary, we make the choice
$ C_{j} = \varepsilon^{-1} 2^{j} $,
$ \varepsilon $ suitable small positive constant. The terms of the form
$ C_{j}^{-1} \| X_{4} \varphi_{N} D_{4}^{N} u\|^{2} $ can be
absorbed on the right hand side of (\ref{SB-E-1.1_D4}).
We have
$\|  \varphi^{(N+1)}_{N}  u \| \leq C^{N +1} \alpha !,$ the analytic growth.
Finally we  observe that the terms in the first sum have the same form as 
$ \| X_{4} \varphi_{N} D_{4}^{N } u \|^{2} $ where one or more 
$ x_{4} $-derivatives have been shifted from $ u $ to $ \varphi_{N} $;
on these terms we can take maximal advantage from the 
sub-elliptic estimate restarting the process. \\
With regard to the second term on the right hand side of  (\ref{es-X_4})
we have 
\begin{equation}\nonumber
\begin{split}
| \langle  x_{1}^{2(r-1)} x_{2}^{2k} \varphi^{(2)}_{N} D_{4}^{N} u,
\varphi_{N} D_{4}^{N} u \rangle |
&\leq 
\frac{1}{2N^{2}} \| X_{4} \varphi^{(2)}_{N} D_{4}^{N -1} u \|^{2} 
+ \frac{N^{2}}{2} \| X_{4} \varphi_{N} D_{4}^{N-1} u \|^{2}   
\\
&\phantom{=}
+| \langle  x_{1}^{r-1} x_{2}^{k} \varphi^{(2)}_{N} D_{4}^{N-1} u,
X_{4} \varphi^{(1)}_{N} D_{4}^{N-1} u \rangle | \nonumber
\\
&\phantom{=}
+ 
| \langle  \frac{1}{N} x_{1}^{r-1} x_{2}^{k} \varphi^{(3)}_{N} D_{4}^{N-1} u,
N X_{4} \varphi_{N} D_{4}^{N-1} u \rangle |\nonumber
\\
&\phantom{=}
+| \langle  x_{1}^{2(r-1)} x_{2}^{2k} \varphi^{(3)}_{N} D_{4}^{N-1} u, 
\varphi^{(1)}_{N} D_{4}^{N-1} u \rangle |. \nonumber
\end{split}
\end{equation}
The last term is the same of the left hand side in which one $ x_{4}$-derivative
has been shifted from  $ u $ to $ \varphi_{N} $ on both side,
we can restart the above process.
On the first two  terms we can take maximal advantage from the 
sub-elliptic estimate restarting the process.
We point out that the `` weight'' $ N $ introduced
above helps to balance the number
of $ x_{4}$-derivatives on $ u $ with the number of derivatives
on $ \varphi_{N} $, we take the factor $ N $ as a derivative on $ \varphi_{N} $
and $ N^{-1} \varphi^{(2)}_{N} $ as $ \varphi^{(1)}_{N}$. 
The other two terms have the same form of the term on the left hand side of
(\ref{es-X_4-1}), the second one with the help of the weight $ N $, we can 
handled both in the same way.    
\vspace{0.2cm}

\noindent
The same strategy can be used to handle the case involving the field $ X_{5} $.

\vspace{0.1cm}

\noindent
The case  $X_{3}$. We have
\begin{equation}\nonumber
 | \langle  \lbrack X_{3}, \varphi_{N} D_{4}^{N} \rbrack u, X_{3}\varphi_{N}
 D_{4}^{N}u \rangle|
\leq  C \| x_{1}^{q-1} \varphi^{(1)}_{N} D_{4}^{N} u \| +
 \frac{1}{C} \| X_{3} \varphi_{N} D_{4}^{N} u\| . 
\end{equation}
The second term  can be absorbed on the left hand side
of (\ref{SB-E-1.1_D4}), if $ C^{-1} $ is chosen small enough.  
Since the first term does not have sufficient power of $ x_{1} $
to take maximal advantage from the sub-elliptic estimate,
we will use the sub-ellipticity. To do this we will pull back $D_{4}^{1/(r+kp)}$.
Let $\chi_{N}(\xi_{4})$ be an Ehrenpreis-H\"ormander cutoff function such 
that $\chi_{N}$ is $ C^{\infty}(\mathbb{R})$ non negative function
such that $\chi_{N} = 0$
for $\xi_{4} < 3 $ and $\chi_{N} =1$ for $\xi_{4} > 4$.
We have
\begin{equation}\nonumber
\| x_{1}^{q-1} \varphi^{(1)}_{N} D_{4}^{N} u \| 
\leq \| x_{1}^{q-1} \varphi^{(1)}_{N} \left( 1-\chi_{N}(N^{-1}D_{4})\right)
D_{4}^{N} u \|
+ \| x_{1}^{q-1} \varphi^{(1)}_{N} \chi_{N}(N^{-1}D_{4})D_{4}^{N} u \|.
\end{equation} 
Since $ 1-\chi_{N}(N^{-1}D_{4})$ has support for $\xi_{4} < 4N$
we can estimate the first therm of the above inequality with
\begin{equation}\nonumber
\| x_{1}^{q-1} \varphi^{(1)}_{N} \left( 1-\chi_{N}(N^{-1}D_{4})\right)
D_{4}^{N} u \|
\leq C^{N+1}N^{N},
\end{equation} 
where $C$ is a positive constant independent by $N$, but depending on $u$.
As already mentioned, to handle the second term of the above inequality
we pull back $D_{4}^{1/(r+kp)}$.
This is well defined since $\xi_{4}> 1$, but is a pseudodifferential operator,
and its commutator with $\varphi_{N}$ needs to some care. 
We use Lemma B.1 and Corollary B.1 in \cite{ABM-TrC1-2016}. 
For completeness we recall them.
Let $\omega_{N} \in C^{\infty}(\mathbb{R})$ be an Ehrenpreis type cutoff 
such that $\omega_{N} =1 $ for $x>2$ and $\omega_{N}=0$ for $x<1$,
$\omega_{N}$ non negative and such that $\omega_{N}\chi_{N}=\chi_{N}$. 
Then we have
\begin{lemma}[\cite{ABM-TrC1-2016}]\label{Lm-Comm}
Let $0<\theta<1$. Then
\begin{equation}\label{Comm-cutoff}
\left[ \omega_{N}\left(N^{-1}D\right)D^{\theta},\varphi_{N}(x)\right]
\chi_{N}\left(N^{-1}D\right)D^{N-\theta}
=\sum_{j=1}^{N} a_{N,j}(x,D)\chi_{N}\left(N^{-1}D\right)D^{N},
\end{equation} 
where $a_{N,j}$ is a pseudo-differential operator of order $-k$ such that
\begin{equation}\label{Es-PDO-Comm}
| \partial_{\xi}^{\alpha}a_{N,k}(x,\xi)|
\leq C_{a}^{j+1} N^{j+\alpha} \xi^{-k-\alpha}, \quad 1\leq j \leq N,
\quad \alpha\leq N.
\end{equation}
\end{lemma}
\begin{corollary}[\cite{ABM-TrC1-2016}]\label{Cr-Comm}
For $1\leq j \leq N-1$ in (\ref{Comm-cutoff}) we have that
\begin{equation}
a_{N,k}(x,D)\chi_{N}\left(N^{-1}D\right)D^{N}
=\frac{\theta(\theta-1) \cdots (\theta-j+1)}{j!} D_{x}^{j}\varphi_{N}(x)
\chi_{N}\left(N^{-1}D\right)D^{N-j}.
\end{equation}
\end{corollary}
Applying these results we find that 
\begin{equation}\nonumber
\begin{split}
\| x_{1}^{q-1} \varphi^{(1)}_{N} \chi_{N}(N^{-1}D_{4})D_{4}^{N} u \|
\leq
&\| x_{1}^{q-1} \varphi^{(1)}_{N} \chi_{N}(N^{-1}D_{4}) 
D_{4}^{N-\frac{1}{r+kp}} u \|_{\frac{1}{r+kp}}\\
&+\sum_{j=1}^{N-1} c_{j} \| x_{1}^{q-1} 
\varphi^{(j+1)}_{N} \chi_{N}(N^{-1}D_{4}) D_{4}^{N-j} u \|\\
&\qquad\qquad\qquad
+ \| x_{1}^{q-1} a_{N,N}(x,D) \chi_{N}(N^{-1}D_{4})D_{4}^{N} u \|.
\end{split}
\end{equation}
The last term has analytic growth. 
To handle the first term on the right hand side we will apply
the sub-elliptic estimate.
Concerning the the terms in the summation, we need, as done previously, 
to pull back $D^{1/(r+kp)}$ once more in order to use the sub-elliptic estimate, 
this will produce either terms with analytic growth or terms of the form 
\begin{equation}\nonumber
c_{j} \| x_{1}^{q-1} \varphi^{(j + 1 )} \chi_{N}(N^{-1}\xi_{4})
D_{4}^{N - j - 1/(r+kp)} u \|_{\frac{1}{r+kp}},
\end{equation}
which can be handled as the first term.

\noindent
Before to analyze the first term on the right hand side of the above inequality
we remark that 
\begin{equation}\nonumber
\begin{split}
| \langle \lbrack \lbrack X_{3}, \varphi_{N} D_{4}^{N} \rbrack ,
X_{3}\rbrack u, \varphi_{N} D_{4}^{N} u \rangle |
&= 
| \langle x_{1}^{2(q-1)} \varphi^{(2)}_{N} D_{4}^{N}  u, 
\varphi_{N} D_{4}^{N}u \rangle
\\
&
 \leq 
\frac{1}{2 N^{2}} \| x_{1}^{q-1} \varphi^{(2)}_{N} D_{4}^{N} u \|^{2} 
+ \frac{N^{2}}{2}\| x_{1}^{q-1} \varphi_{N} D_{4}^{N} u \|^{2}.
\end{split}
\end{equation}
As above we use the ``weight '' $ N $ to balance the number
of $ x_{4}-$ de\-ri\-va\-ti\-ves on $ u $ with the number of
derivatives on $ \varphi_{N} $. 
The two terms on the right hand side have the same form
as $ \|  x_{1}^{q-1} \varphi^{(1)}_{N} D_{4}^{N} u \| $,
we can use the same strategy 
to analyze these two terms.\\
Then the only term that we have to handle is the term 
$ \| x_{1}^{q-1} \varphi^{(1)}_{N} \chi_{N}(N^{-1}\xi_{4})
D_{4}^{N -\frac{1}{r+pk}}u\|_{1/r+pk}$.
Once again, to estimate this term we use the sub-elliptic estimate
(\ref{SB-E-1/2}) 
replacing $ u $ with $ x_{1}^{q-1} \varphi^{(1)}_{N}\chi_{N}(N^{-1}\xi_{4})
D_{4}^{N -\frac{1}{r+pk}} u $. We have
\begin{equation} \nonumber
\begin{split}
&\|  x_{1}^{q-1}\varphi^{(1)}_{N} \chi_{N}(N^{-1}\xi_{4})
D_{4}^{N-\frac{1}{r+kp}} u \|^{2}_{\frac{1}{r+kp}}
\!\!\!+ \!\!\sum_{j =1}^{5} \|X_{j}  x_{1}^{q-1}\varphi^{(1)}_{N}
\chi_{N}(N^{-1}\xi_{4}) D_{4}^{N-\frac{1}{r+kp}} u \|^{2}
\\
&\quad\leq
 \| x_{1}^{q-1}\varphi^{(1)}_{N} \chi_{N}(N^{-1}\xi_{4})
 D_{4}^{N-\frac{1}{r+kp}} P u \|^{2}
+ \| x_{1}^{q-1}\varphi^{(1)}_{N} \chi_{N}(N^{-1}\xi_{4})
D_{4}^{N-\frac{1}{r+kp}} u\|^{2}
\\
&\qquad
+2 \sum_{j=1}^{5} |\langle X_{j}[X_{j}, x_{1}^{q-1}\varphi^{(1)}_{N}]
\chi_{N}(N^{-1}\xi_{4}) D_{4}^{N-\frac{1}{r+kp}} u , 
 x_{1}^{q-1}\varphi^{(1)}_{N} \chi_{N}(N^{-1}\xi_{4})
 D_{4}^{N-\frac{1}{r+kp}} u\rangle|
\\
&\qquad
+  \sum_{j=1}^{5} |\langle [X_{j},[X_{j}, x_{1}^{q-1}\varphi^{(1)}_{N}] ]
\chi_{N}(N^{-1}\xi_{4}) D_{4}^{N-\frac{1}{r+kp}} u , 
 x_{1}^{q-1}\varphi^{(1)}_{N}\chi_{N}(N^{-1}\xi_{4})
 D_{4}^{N-\frac{1}{r+kp}} u\rangle|.
\end{split}
\end{equation}
The right hand side of the above equation can be estimate by
\begin{equation}\label{est3}
\begin{split}
& 
C(q-1)^{2}\| x_{1}^{q-2} \varphi^{(1)}_{N} \chi_{N}(N^{-1}\xi_{4})
D_{4}^{N -\frac{1}{r+kp}} u \|^{2}
+2C \| x_{1}^{2(q-1)} \varphi^{(2)}_{N}\chi_{N}(N^{-1}\xi_{4})
D_{4}^{N -\frac{1}{r+kp}} u \|^{2}
\\
&
+\frac{1}{N^{2}} \| x_{1}^{2(q-1)} \varphi^{(3)}_{N}\chi_{N}(N^{-1}\xi_{4})
D_{4}^{N-\frac{1}{r+kp}} u \|^{2}
+N^{2}\| x_{1}^{2(q-1)} \varphi^{(1)}_{N} \chi_{N}(N^{-1}\xi_{4})
D_{4}^{N -\frac{1}{r+kp}} u \|^{2}
\\
&
+2C \| x_{1}^{r+q-2} x_{2}^{k}\varphi^{(2)}_{N} \chi_{N}(N^{-1}\xi_{4})
D_{4}^{N -\frac{1}{r+kp}} u \|^{2}
\\
&
+2C \| x_{1}^{r+\ell+q-2} \varphi^{(2)}_{N}\chi_{N}(N^{-1}\xi_{4})
D_{4}^{N -\frac{1}{r+kp}} u \|^{2}
\\
&
+ | \langle x_{1}^{2(r-1)+q-1} x_{2}^{2k}
\varphi^{(3)}_{N}\chi_{N}(N^{-1}\xi_{4}) 
D_{4}^{N -\frac{1}{r+kp}} u,
x_{1}^{q-1} \varphi^{(1)}_{N}\chi_{N}(N^{-1}\xi_{4})
D_{4}^{N -\frac{1}{r+kp}} u \rangle 
\\
&
+  | \langle x_{1}^{2(r+\ell-1)+q-1}\varphi^{(3)}_{N}\chi_{N}(N^{-1}\xi_{4})
D_{4}^{N -\frac{1}{r+kp}} u,
x_{1}^{q-1} \varphi^{(1)}_{N} \chi_{N}(N^{-1}\xi_{4})
D_{4}^{N -\frac{1}{r+kp}} u \rangle
\end{split}
\end{equation}
modulo terms which can be absorbed on the left hand side
or which give analytic growth.
We remark that on the last four terms we can take maximal advantage from the
sub-elliptic estimate restarting the processes; moreover in view of the role of 
the weight $ N $ the third and the fourth term have the same form of the 
second one.    
Before to give the general form of the terms which appear inside of the iterating
process we analyze the particular situations. 
To handle the first term on the right hand side of (\ref{est3}) 
we must use the sub-ellipticity, i.e. we pull back $D_{4}^{1/(r+kp)}$.  
Using the Lemma \ref{Lm-Comm} and the Corollary \ref{Cr-Comm} we have
\begin{equation}\nonumber
\begin{split} 
\| x_{1}^{q-2} \varphi^{(1)}_{N} \chi_{N}(N^{-1}D_{4})
D_{4}^{N-\frac{1}{r+kp}} u \|
\leq
&\| x_{1}^{q-2} \varphi^{(1)}_{N} \chi_{N}(N^{-1}D_{4})
D_{4}^{N-\frac{2}{r+kp}} u \|_{\frac{1}{r+kp}}\\
&+\sum_{j=1}^{N-1} c_{j} \| x_{1}^{q-2} \varphi^{(j+1)}_{N} \chi_{N}(N^{-1}
D_{4}) D_{4}^{N-j-\frac{1}{r+kp}} u \|\\
&\qquad\qquad\quad+ \| x_{1}^{q-2} a_{N,N}(x,D) \chi_{N}(N^{-1}
D_{4})D_{4}^{N-\frac{1}{r+kp}} u \|.
\end{split}
\end{equation}
The last term has analytic growth. 
To handle the first term on the right hand side 
we will apply the sub-elliptic estimate.
Concerning the the terms in the summation, we need, as done previously, 
to pull back $D^{1/(r+kp)}$ once more in order to use the sub-elliptic estimate, 
this will produce either terms with analytic growth or terms of the form 
\begin{equation}\nonumber
c_{j} \| x_{1}^{q-1} \varphi^{(j+1)} \chi_{N}(N^{-1}D_{4})
D_{4}^{N - j - \frac{2}{(r+kp)}} u \|_{\frac{1}{r+kp}},
\end{equation}
which can be handled as the first term.\\
\noindent
Iterating the above strategy at the $j$-th step we obtain a term of the form
\begin{equation}\nonumber
\| x_{1}^{q-j-1} \varphi^{(1)}_{N} \chi_{N}(N^{-1}\xi_{4})
D_{4}^{N -\frac{j+1}{r+kp}} u \|^{2}_{\frac{1}{r+kp}}.
\end{equation}
When $ j = q-1$ we have $ \| \varphi^{(1)}_{N} \chi_{N}(N^{-1}\xi_{4})
D_{4}^{N -\frac{q}{r+kp}} u \|^{2}_{1/r+kp}$.
Iterating this cycle $s$-times we obtain a term of the form
\begin{equation}\nonumber
C^{s} \| \varphi^{(s)}_{N} \chi_{N}(N^{-1}\xi_{4})D_{4}^{N -s\frac{q}{r+kp}}
u \|^{2}_{\frac{1}{r+kp}}.
\end{equation}
Using up all $x_{4}$-derivatives we estimate this term,
hence the right hand side
of (\ref{SB-E-1.1_D4}), with $ C^{2(N +1)} N^{2N(r+kp)/q}$.
We have a growth corresponding to $ G^{(r+kp)/q}$.\\
\noindent 
The second term on the right hand side of (\ref{est3}),
$\| x_{1}^{2(q-1)} \varphi^{(2)}_{N}\chi_{N}(N^{-1}\xi_{4})
D_{4}^{N -\frac{1}{r+kp}} u \|^{2}$, once again 
we must use the sub-ellipticity, that is using the Lemma \ref{Lm-Comm} and
the Corollary \ref{Cr-Comm} we pull back $D_{4}^{1/(r+kp)}$
restarting the process.\\ 
Iterating this strategy at the $ h$-th step we obtain a term of the form
\begin{equation}\nonumber
\| x_{1}^{h(q-1)} \varphi^{(h)}_{N}\chi_{N}(N^{-1}\xi_{4})
D_{4}^{N -\frac{h-1}{r+kp}} u \|^{2}.
\end{equation}
Let $ 0<\beta<1 $ a  parameter that will be chosen later.
Using the Lemma \ref{Lm-Comm} and the Corollary \ref{Cr-Comm} we pull back 
$D_{4}^{\beta}$; we can estimate the above quantity with
\begin{equation}\label{Est-beta}
\| x_{1}^{h(q-1)}D_{4}^{\beta} \varphi^{(h)}_{N}\chi_{N}(N^{-1}\xi_{4})
D_{4}^{N-\frac{h-1}{r+kp}-\beta} u \|^{2}
\end{equation}
modulo terms of the form 
$c_{j}\| x_{1}^{h(q-1)}D_{4} \varphi^{(h+j)}_{N}\chi_{N}(N^{-1}\xi_{4})
D_{4}^{N-\frac{h-1}{r+kp}-j} u \|^{2}$, 
$j = 1,\dots, N-1$, and $\| x_{1}^{h(q-1)}a_{N,N}(x,D)\chi_{N}(N^{-1}\xi_{4})
D_{4}^{N-\frac{h-1}{r+kp}} u \|^{2}$. The last one gives analytic growth, the 
others can be estimated restarting the process,
i.e. pulling back $D_{4}^{\beta}$
and using the same process to estimate (\ref{Est-beta}),
that we will show below.
The term (\ref{Est-beta}) can be estimated by
\begin{equation}\nonumber
\begin{split}
&\| x_{1}^{h(q-1)-(p-1)} x_{2} D_{4}^{2\beta} \varphi^{(h)}_{N}
\chi_{N}(N^{-1}\xi_{4})
D_{4}^{N -\frac{h-1}{r+kp}-\beta} u \|^{2}
+ \| X_{2} \varphi^{(h)}_{N}\chi_{N}(N^{-1}\xi_{4})
D_{4}^{N -\frac{h-1}{r+kp}-\beta} u \|^{2}
\\
&\qquad\quad
\leq
C_{1} \| x_{1}^{(m+1)h(q-1)-m(p-1)} x_{2}^{m} D_{4}^{(m+1)\beta}
\varphi^{(h)}_{N}\chi_{N}(N^{-1}\xi_{4})
D_{4}^{N-\frac{h-1}{r+kp}-\beta} u \|^{2}
\\
&\qquad\qquad
+ \| X_{2} \varphi^{(h)}_{N}\chi_{N}(N^{-1}\xi_{4})
D_{4}^{N-\frac{h-1}{r+kp}-\beta} u \|^{2}
+C_{2} \| x_{1}^{h(q-1)} D_{4}^{\beta} \varphi^{(h)}_{N}\chi_{N}
(N^{-1}\xi_{4})D_{4}^{N-\frac{h-1}{r+kp}-\beta} u \|^{2}
\end{split}
\end{equation}
where $ C_{2} $ is a small suitable constant. The last term can be absorbed 
on the left hand side. Choosing $ m = k $, $ \beta = ( k+1 )^{-1} $ and
$  h=( r -1 + k( p-1))/ ( (k+1) (q-1)) $ we obtain 
\begin{equation}\nonumber
\begin{split}
\| X_{4} \varphi^{(h)}_{N}\chi_{N}(N^{-1}\xi_{4})
D_{4}^{N -\frac{h-1}{r+kp}-\beta} u \|^{2}
+ \| X_{2} \varphi^{(h)}_{N}\chi_{N}(N^{-1}\xi_{4})
D_{4}^{N-\frac{h-1}{r+kp}-\beta} u \|^{2}.
\end{split}
\end{equation}
Restarting the process, taking maximum advantage from the sub-elliptic
estimate we obtain after $ s $ step
\begin{equation}\nonumber
\begin{split}
\| X_{4} \varphi^{(sh)}_{N}\chi_{N}(N^{-1}\xi_{4})
D_{4}^{N -s\left(\frac{h-1}{r+kp}+\beta\right)} u \|^{2}
+ \| X_{2} \varphi^{(sh)}_{N}\chi_{N}(N^{-1}\xi_{4})
D_{4}^{N-s\left(\frac{h-1}{r+kp}+\beta\right)} u \|^{2}.
\end{split}
\end{equation}
Iterating until all the $ x_{4}$-derivatives are used up,
that is until $N -s((h-1)(r+kp)^{-1}+\beta) \sim 0 $,
we have the growth corresponding to $ G^{(r+kp)/q}$.\\
Combining and iterating the above processes more time,
removing powers of $ x_{1}$ and $ x_{2}$ with $ D_{4}$
and taking profit from the sub-ellipticity we may estimate
the left hand side of (\ref{SB-E-1.1_D4}) with terms of the form
\begin{equation}\nonumber
\begin{split}
&N^{-2m_{0}}
\Big(
 \| X_{4}\varphi^{(m_{1}+hm_{2}+m_{3}+m_{0})}_{N} 
\chi_{N}(N^{-1}\xi_{4})
D_{4}^{N -(m_{1}+m_{2}h)\frac{q}{r+pk}-m_{3} }u \|^{2}\\
&
+\| X_{2}\varphi^{(m_{1}+hm_{2}+m_{3}+m_{0})}_{N} 
\chi_{N}(N^{-1}\xi_{4})
D_{4}^{N -(m_{1}+m_{2}h)\frac{q}{r+pk}-m_{3} }u \|^{2}\\
&
\left.
+\| x_{1}^{m_{4}(q-1) -m_{5}}
\varphi^{(m_{0}+ m_{1}+ m_{2}h + m_{3}+m_{4})}_{N} 
\chi_{N}(N^{-1}\xi_{4})
D_{4}^{N - m_{1}+\frac{(m_{3}+m_{2}h)q}{r+kp}-\frac{m_{4} + m_{5}}{r+kp}
- m_{6} +\frac{m_{6}}{r+kp} }u\|^{2}
u \|^{2}_{\frac{1}{r+kp}}
\right)
\end{split}
\end{equation}
where $ h$ is as above, $ (h-1)(r+kp)^{-1}-\beta= qh(r+kp)^{-1}$
and $ m_{4}(q-1) -m_{5} \leq (q-1)$.
Iterating until all $ x_{4}$-derivatives are used up, that is
$N- (m_{1}+m_{2}h)q(r+pk)^{-1}-m_{3}\sim 0 $ and
$ N - m_{1}+(m_{3}+m_{2}h)q(r+kp)^{1}- 
(m_{4} + m_{5})(r+kp)^{1} - m_{6} +m_{6}(r+kp)^{-1}  \sim 0$
we have that
$ m_{1}+hm_{2}+m_{3}$ and $ m_{1}+hm_{2}+m_{3}+m_{4}$ ,
since $m_{3} \geq 1 $ and $ m_{6} \geq 1$,
are small or equal to 
$ (r+kp)N q^{-1}$. We can conclude
\begin{equation}\nonumber
\| \varphi_{N} D^{N}_{4} u\|^{2}_{\frac{1}{r+kp}}\! 
+\! \sum_{j=1}^{5} \|X_{j} \varphi_{N} D^{N}_{4} u \|^{2}
\leq 
C^{2(N+1)} (N)^{2N\frac{r+kp}{q}}
\end{equation}
where $ C$ is independent by $N$ but depends on $ u$.
This conclude the proof.
\begin{remark}
In particular we have that if $u$ solves the equation
$P_{1}u = f$ and $f$ is analytic, if $\rho_{_{0}}\in \Sigma_{r+kp}$
then  $\rho_{_{0}}\notin WF_{\frac{r+kp}{q}}(u)$ and if
$\rho_{_{1}}\in\Sigma_{r}$ then $\rho_{_{1}}\notin WF_{\frac{r}{p}}(u)$.
\end{remark}
\subsection{Gevrey Regularity for $P_{2}(x;D)$}
\noindent
The characteristic variety of $P_{2}$ is:
\begin{equation}\nonumber
Char(P_{2}) 
=\lbrace (x;\xi)\in T^{*}\mathbb{R}^{4}\setminus \lbrace 0 \rbrace\,
:\, x_{1}=0,\, \xi_{1}=0\rbrace.
\end{equation}
\noindent
Following the ideas in \cite{Tr-2006} and \cite{BoTr} it can be seen as the 
disjoint union of analytic submanifolds: 
\begin{equation}\nonumber
Char(P_{2}) = \Sigma_{p,\pm} \cup \Sigma_{q,\pm} \cup \Sigma_{r,\pm}
\cup \Sigma_{r+qk,\pm},
\end{equation}
where
\begin{equation}\nonumber
\begin{array}{rcl}
\Sigma_{p,\pm}&=&
\lbrace  (x;\xi) \in T^{*}\mathbb{R}^{4}\setminus \lbrace 0
\rbrace :\, x_{1}=0,\, \xi_{1}=0, \,\xi_{2}\gtrless 0 \rbrace,
\\[5pt]
\Sigma_{q,\pm}  &=& 
\lbrace  (x;\xi)\in T^{*}\mathbb{R}^{4}\setminus \lbrace 0
\rbrace :\, x_{1}=0,\, \xi_{1}=0,\, \xi_{2} =0, \,\xi_{3}\gtrless 0 \rbrace,
\\[5pt]
\Sigma_{r,\pm} &=& 
\lbrace  (x;\xi)\in T^{*}\mathbb{R}^{4}\setminus \lbrace 0
\rbrace :\, x_{1}=0,\, \xi_{1}=0,\, \xi_{2} =0,\, \xi_{3}=0,
\,\xi_{4}x_{3}\neq 0 
\rbrace
\\[5pt]
\Sigma_{r+qk,\pm} &=&  
\lbrace  (x;\xi)\in T^{*}\mathbb{R}^{4}\setminus \lbrace 0 
\rbrace :\,
x_{1}=0,\, \xi_{1}=0,\, \xi_{2} =0,\, \xi_{3}=0,
\, x_{3}= 0,\, \xi_{4}\gtrless 0 \rbrace.
\end{array}
\end{equation}
\noindent
Once more by the results in \cite{DZ} and \cite{RS} the operator $ P_{2} $ has the following sub-
elliptic estimate with loss of $2(1-1/(r+kq))$ derivatives:
%
\begin{equation}
\label{SB-E-1/3}
\|u\|^{2}_{\frac{1}{r+kq}} + \sum_{j=1}^{5} \|X_{j} u\|^{2} 
\leq 
C \big( |\langle P_{2}u, u \rangle | + \| u \|^{2} \big). 
\end{equation}
Here
$ X_{1} = D_{1} $, $ X_{2} = x_{1}^{p-1}D_{2}$,
$ X_{3} = x_{1}^{q-1} D_{3}$,
$ X_{4} = x_{1}^{r-1} x_{3}^{k} D_{4} $, $ X_{5} = x_{1}^{r+\ell-1} D_{4}$,
$ \| \cdot \|_{s} $ denotes the $ H^{s} $ Sobolev norm and
$ \| \cdot \| = \| \cdot \|_{0} $
denotes the $ L^{2} $ norm on the fixed open set $ \Omega $.\\   
To study the regularity of the solutions we estimate the high order derivatives 
of the solutions in $ L^{2} $ norm, as in the case of $P_{1}$.
For $ x_{1} \neq 0 $ the operator $ P_{2} $ is elliptic and we shall not examine
this region, elliptic operators are Gevrey hypoelliptic in any class $ G^{s} $ for
$ s \geq 1 $.\\
Let $ \varphi_{N} (x) $ be a cutoff function of Ehrenpreis-H\"or\-man\-der
type with the same properties described in the beginning of the previous 
paragraph.  

We assume that $ \varphi_{N} $ is independent of  the $ x_{1} $-variable
for the same reason described in the proof of the regularity of $P_{1}$.
Moreover we may assume that $ \varphi_{N} $
is independent of the $ x_{3} $-variable since every $ x_{3}$-derivative
landing on $ \varphi_{N} $ would leave a cut off function supported where
$ x_{3} $ is bounded away from zero, in this region the operator 
satisfies the H\"or\-man\-der-Lie algebra condition at the step $ r $.
The operator $ P_{2} $ is sub-elliptic with loss of $ 2 ( 1 - 1/r ) $ derivatives.
In this region the operator is a generalization of  the Ole\u{\i}nik-Radkevi\v c
operator then $ P_{2}$ is $ G^{r/p}$-hypoelliptic and not better,
for more details see \cite{C} and \cite{BT}.
Thus we can conclude that 
if $v$ solves the equation $P_{1}v = f$ and $f$ is analytic then the points
$\displaystyle\rho_{_{1}}= \in \Sigma_{r}$, does not belong to
$WF_{\frac{r}{p}}(v)$.\\
Now, we are interested to the microlocal regularity 
in  $\displaystyle\rho_{_{0}}\in \Sigma_{r+kp}$, to obtain
this it is sufficient to study the microlocal regularity of $P_{2}$
in $(0;e_{4})$.
The microlocal regularity in  a generic point $\displaystyle\rho_{_{0}}$ can be 
obtained following the same strategy below  with the only difference that the 
cut-off function $\varphi_{N}(x)$ will be identically equal to 1 in $\Omega_{0}$ 
neighborhood of $\pi_{x}(\rho_{_{0}})=(0,x_{2}^{0},0,x_{4}^{0})$.
Thus since we are interested to the microlocal regularity of $P_{2}$
in $(0;e_{4})$ we take
$\varphi_{N} (x) = \varphi_{N}(x_{2},x_{4})$.\\
%
%
We replace $ u $ by $ \varphi_{N} (x) D_{4}^{N} u $ in ( \ref{SB-E-1/3}).
We have
\begin{equation}\label{SB-E-1.2_D4}
\| \varphi_{N} D^{N}_{4} u\|^{2}_{\frac{1}{r+kq}}\!
+\! \sum_{j=1}^{5} \|X_{j} 
\varphi_{N} D^{N}_{4} u \|^{2}
\leq \!
C \left( |\langle P_{2} \varphi_{N} D^{N}_{4} u, \varphi_{N} D^{N}_{4} u 
\rangle | + \| \varphi_{N} D^{N}_{4} u \|^{2}  \right). 
\end{equation}
As in the case of the operator $P_{1}$ we want to estimate terms of the form:
\begin{equation}
\begin{split}
 \langle  \lbrack X_{j}, \varphi_{N} D_{4}^{N} \rbrack u,
 X_{j}\varphi_{N} D_{4}^{N} u \rangle
\text{ and }
 \langle \lbrack \lbrack X_{j}, \varphi_{N} D_{4}^{N} \rbrack , X_{j}\rbrack u, 
 \varphi_{N} D_{4}^{N} u \rangle,
\qquad j=1,2,3,4,5.  
\end{split}
\end{equation}
Since $ \varphi_{N} $ depends only by $ x_{2} $ and $ x_{4} $
we must analyze the commutators with $ X_{2}$, $ X_{4}$ and $ X_{5}$.
The cases $ X_{4} $ and $X_{5}$ give analytic growth, they can handled in 
same way as done in the study of $P_{1}$;
in these cases we can take maximal advantage from the sub-elliptic estimate.
The case $X_{2}$. In this case we have to estimate the term
\begin{equation}\nonumber
 \| x_{1}^{p-1} \varphi^{(1)}_{N} D_{4}^{N} u \|. 
\end{equation}
Since it does not have sufficient power of $ x_{1} $
to take maximal advantage from the sub-elliptic estimate,
we will use the sub-ellipticity. To do this we will pull back $D_{4}^{1/(r+kq)}$.
Using the same strategy employed to study the case 
of the vector field $X_{3}$
in the study of the regularity of $P_{1}$, here we have $x_{1}^{p-1}$
instead of $x_{1}^{q-1}$.
Following the same strategy used to deduce the regularity of $P_{1}$,
we  conclude that
\begin{equation}\nonumber
\| \varphi_{N} D^{N}_{4} u\|^{2}_{\frac{1}{r+kq}}\! 
+\! \sum_{j=1}^{5} \|X_{j} 
\varphi_{N} D^{N}_{4} u \|^{2}
\leq 
C^{2(N+1)} (N)^{2N\frac{r+kq}{p}},
\end{equation}
where $ C$ is independent by $N$ but depends on $ u$.
We have that the point 
$(0;e_{4})$ and more in general that the points $\rho_{0} \in 
\Sigma_{r+kq}$ do not belong to $WF_{\frac{r+kq}{p}}(u)$.
This conclude the proof of the theorem.
\begin{remark}
In particular we have that if $u$ solves the equation
$P_{1}u = f$ and $f$ is analytic, if $\rho_{_{0}}\in \Sigma_{r+kp}$
then  $\rho_{_{0}}\notin WF_{\frac{r+kp}{p}}(u)$ and if
$\rho_{_{1}}\in\Sigma_{r}$ then $\rho_{_{1}}\notin WF_{\frac{r}{p}}(u)$.
\end{remark}
\section{Proof of Theorem \ref{T-2}}\label{P-T_2}
The characteristic variety of $P_{3}$ is:
\begin{equation}\nonumber
Char(P_{3}) =\lbrace (x;\xi)\in T^{*}\mathbb{R}^{4}\setminus \lbrace 0 
\rbrace\, :\, x_{1}=0,\, \xi_{1}=0\rbrace.
\end{equation}
\noindent
Following the ideas in \cite{Tr-2006} and \cite{BoTr} it can be seen as the 
disjoint union of analytic submanifolds: 
\begin{equation}\nonumber
\begin{array}{rcl}
\Sigma_{p,\pm}&=&
\lbrace  (x;\xi) \in T^{*}\mathbb{R}^{4}\setminus \lbrace 0
\rbrace :\, x_{1}=0,\, \xi_{1}=0, \,\xi_{2}\gtrless 0 \rbrace,
\\[5pt]
\Sigma_{q,\pm}  &=& 
\lbrace  (x;\xi)\in T^{*}\mathbb{R}^{4}\setminus \lbrace 0
\rbrace :\, x_{1}=0,\, \xi_{1}=0,\, \xi_{2} =0, \,\xi_{3}\gtrless 0 \rbrace,
\\[5pt]
\Sigma_{r,\pm} &=& 
\lbrace  (x;\xi)\in T^{*}\mathbb{R}^{4}\setminus \lbrace 0
\rbrace :\, x_{1}=0,\, \xi_{1}=0,\, \xi_{2} =0,\, \xi_{3}=0,
\,\xi_{4}x_{2}\neq 0 
\rbrace,
\end{array}
\end{equation}
if $f >pk+r$ there is only one more stratum of depth $r+kp$:
\begin{equation}\nonumber
\Sigma_{r+pk,\pm} =  
\lbrace  (x;\xi)\in T^{*}\mathbb{R}^{4}\setminus \lbrace 0 
\rbrace :\, 
x_{1}=0,\,\xi_{1}=0,\, \xi_{2} =0,\,\, x_{2}= 0, \xi_{3}=0,
\, \xi_{4}\gtrless 0 \rbrace ;
\end{equation}
otherwise if $f <pk+r$ there are other two strata of depth $f$,
\begin{equation}\nonumber
\begin{array}{rcl}
\Sigma_{f,\pm} \!&=&\!   
\lbrace  (x;\xi) \in T^{*}\mathbb{R}^{4}\setminus \lbrace 0
\rbrace :\, x_{1}=0,\, \xi_{1}=0,\, \xi_{2} =0,\,\, x_{2}= 0, \xi_{3}=0,\, x_{3}
\xi_{4}\neq 0\rbrace,
\end{array}
\end{equation}
and of depth $r+kp$ if $r+kp < f+q\ell$, 
\begin{equation}\nonumber
\begin{array}{rcl}
\Sigma_{r+kp,\pm}
\!&=&\! 
\lbrace  (x;\xi)\in T^{*}\mathbb{R}^{4}\setminus \lbrace 0
\rbrace :\, x_{1}=0,\, \xi_{1}=0,\, \xi_{2} =0,\,\, x_{2}= 0, \xi_{3}=0,\, x_{3} 
=0,\, 
\xi_{4}\gtrless 0 \rbrace,
\end{array}
\end{equation}
or of depth $f + q\ell$ if $r+kp > f+q\ell$,
\begin{equation}\nonumber
\begin{array}{rcl}
\Sigma_{f+q\ell,\pm}
\!&=&\! 
\lbrace  (x;\xi)\in T^{*}\mathbb{R}^{4}\setminus \lbrace 0
\rbrace :\, x_{1}=0,\, \xi_{1}=0,\, \xi_{2} =0,\,\, x_{2}= 0, \xi_{3}=0,\, x_{3} 
=0,\, 
\xi_{4}\gtrless 0 \rbrace.
\end{array}
\end{equation}
\noindent
{\bf  Case $ \textit{f} > \textit{r}+\textit{kp}$}:
In this case the H\"ormander condition is satisfied at the step $r+kp$.
Once more by the results in \cite{DZ} and \cite{RS}
$ P_{3}$ satisfies the following sub-elliptic estimate 
with loss of $2(1-1/(r+kp))$ derivatives:
\begin{equation}\label{SE-E-P-3}
\|u\|^{2}_{\frac{1}{r+kp}} + \sum_{j=1}^{6} \|X_{j} u\|^{2} 
\leq 
C \big( |\langle P_{3}u, u \rangle | + \| u \|^{2} \big). 
\end{equation}
Here $ X_{1} =\! D_{1} $, $ X_{2} = \!x_{1}^{p-1}D_{2}$,
$ X_{3} =\! x_{1}^{q-1} D_{3}$,
$ X_{4} = \!x_{1}^{r-1} x_{2}^{k} D_{4} $,
$ X_{5} =\! x_{1}^{f-1}x_{3}^{\ell} D_{4}$
and $X_{6} =\! x_{1}^{f+e-1}D_{4}$.\\
The result can be archived following the some strategy used to characterize the 
regularity of the operator $P_{1}(x,D)$, Theorem \ref{T-1}--\textbf{\textit{i}}.
In fact the presence of the additional
vector field $ X_{5} =\! x_{1}^{f-1}x_{3}^{\ell} D_{4}$ gives, in the algorithm 
developed to handle the operator $P_{1}$, only a negligible contribution,
i.e. analytic growth: to estimate the terms 
$| \langle x_{1}^{f-1}x_{3}^{\ell} \varphi^{(1)}_{N} D_{4}^{ N } u,
X_{5}\varphi_{N} D_{4}^{N}u \rangle |$
and
$| \langle  x_{1}^{2(f-1)} x_{3}^{2\ell} \varphi^{(2)}_{N} D_{4}^{N} u,
\varphi_{N} D_{4}^{N} u \rangle | $
can take maximal advantage from the sub-elliptic estimate. 

\vspace{0.3cm}

\noindent
{\bf Case $ \textit{f} < \textit{r}+\textit{kp}$}:
In this case we distinguish two different situations: $ r+kp < f +\ell q$
and $ f +\ell q <  r+kp$.
Since the only difference between the two cases is the subelliptic index, that is 
in the first case the H\"ormander condition is satisfied at the step $r+kp$ and in 
the other at the step $f +\ell q$ we will analyze only the first one.\\
{\bf Case $ \textit{r}+\textit{kp} < \textit{f} +\mathbf{\ell}\textit{q}\,$}: 
The operator $P_{3}$ is
sub-elliptic with loss of $2(1-1/(r+kp))$ derivatives, as
above the sub-elliptic a priori estimate (\ref{SE-E-P-3}) holds.\\
Let $\varphi_{N}(x)$ be a localizing 
cutoff function of Ehrenpreis-H\"ormander type.
We may assume that $\varphi_{N}$ is independent of the  $x_{1}$-variable
since every $ x_{1} $-derivative landing on $ \varphi_{N} $ would leave a
cutoff function supported where $ x_{1} $ is bounded away from zero, 
where the operator is elliptic. We can also assume that $\varphi_{N}$ is 
independent of the  $x_{2}$-variable. 
If $x_{2}\neq 0$ the operator $P_{3}$ is an operator
of Ole\u{\i}nik-Radkevi\v c type, \cite{OR1973},
in view of the result obtained in \citep{C}, in this region, the operator is
$G^{r/p}$-hypoelliptic. 
We can conclude that if $u$ solves the equation $P_{3}u = g$ 
and $g$ is analytic then the points $\rho_{_{6}} \in \Sigma_{r}$, 
does not belong to $WF_{r/p}(u)$.\\ 
Moreover we may assume that $\varphi_{N}$ is independent of the
$x_{3}$-variable.
Every $x_{3}$-derivative landing on $\varphi_{N}$ would leave a cut off 
function supported where $x_{3}$ is bounded away from zero,
in this region the H\"ormander condition is satisfied at the step $f$.
The operator $P_{3}$ has the same form of the operator $P_{1}$, 
(\ref{Op_1}), in the Theorem \ref{T-1}, with $pk>\ell$.
We can conclude that if $u$ solves the equation $P_{3}u = g$ and $g$
is analytic then the points $\rho_{_{5}} \in \Sigma_{f}$,
do not belong to $WF_{f/p}(u)$.\\
We assume that $\varphi_{N}(x)= \varphi_{N}(x_{4})$.
We replay $u$ by $\varphi_{N}(x_{4})D_{4}^{N}u$ in (\ref{SE-E-P-3}).
We have 
\begin{equation}\label{RS-P_1-ii-D4}
\| \varphi_{N} D^{N}_{4} u\|^{2}_{\frac{1}{r+kp}}\! +
\! \sum_{j=1}^{6} \|X_{j} \varphi_{N} D^{N}_{4} u \|^{2}
\leq \! C 
\left( |\langle P_{3} \varphi_{N} D^{N}_{4} u, \varphi_{N} D^{N}_{4} u\rangle | 
+ \| \varphi_{N} D^{N}_{4} u \|^{2}  \right). 
\end{equation}
We have to estimate terms of the form:
\begin{equation}\nonumber
 \langle  \lbrack X_{j}, \varphi_{N} D_{4}^{N} \rbrack u, X_{j}
 \varphi_{N} D_{4}^{N}u \rangle
\text{ and }
 \langle \lbrack \lbrack X_{j}, \varphi_{N} D_{4}^{N} \rbrack , X_{j}\rbrack u, 
 \varphi_{N} D_{4}^{N} u \rangle, \qquad j=1,2,3,4,5, 6.  
\end{equation}
Since $ \varphi_{N} $ depends only by $ x_{4} $, $ X_{1}$, $ X_{2}$ and
$ X_{3}$ commute with $\varphi_{N}$. We must only analyze the commutators 
with $ X_{4}$, $ X_{5}$ and $ X_{6}$.
These cases give analytic growth,
we can take maximal advantage from the sub-elliptic estimate. 
They can be handled as the field $X_{4}$, (\ref{es-X_4}),
in the proof of the Theorem \ref{T-1}.
We conclude that the point $(0;e_{4})$ and more in general the points 
$\rho_{_{4}}\in \Sigma_{r+kp}$ do not belong to $WF_{a}(u)$.
\begin{remark}
In particular we have that if $u$ solves the equation
$P_{3}u = f$ and $f$ is analytic, if $f > r + kp$
then if $\displaystyle\rho_{_{2}}\in \displaystyle\Sigma_{r+kp}$
then $\displaystyle\rho_{_{2}}\notin WF_{\frac{r+kp}{q}}(u)$ and 
if $\displaystyle\rho_{_{3}}\in \Sigma_{r}$ then
$ \displaystyle\rho_{_{3}}\notin WF_{\frac{r}{p}}(u)$;
if $f \leq r + kp$ then if
$\displaystyle\rho_{_{4}}\in \displaystyle\Sigma_{r+kp}$
then $\displaystyle\rho_{_{4}}\notin WF_{a}(u)$, 
if $\displaystyle\rho_{_{5}}\in \Sigma_{f}$
then $\displaystyle\rho_{_{5}}\notin WF_{\frac{f}{p}}(u)$ and 
if $\displaystyle\rho_{_{6}}\in \Sigma_{r}$
then $\displaystyle\rho_{_{6}}\notin WF_{\frac{r}{p}}(u)$.            
\end{remark}
%
%
\section{On the partial regularity of the operators $\textit{\textbf{P}}_{\!\text{\textbf{1}}}$ and $\textit{\textbf{P}}_{\!\text{\textbf{2}}}$}\label{PR_P_1&P_2}
In this section, following the ideas in \cite{BT}, we analyze the partial regularity in a neighborhood of the origin 
of the operators $P_{1}$, (\ref{Op_1}), and $P_{2}$, (\ref{Op_2}). 
We recall the definition of the non-isotropic Gevrey classes:
\begin{definition}
A smooth function $f(x_{0},x_{1},\dots,x_{n})$ belongs to the Gevrey space $G^{(\alpha_{0},\alpha_{1},\dots,\alpha_{n})}$
at the point $x_{0}$ provided that there exists a neighborhood, $ U $,
of $x_{0}$ and a constant $C_{f}$ such that for all multi-indices $\beta$ 
$$
|D^{\beta}f| \leq C_{f}^{|\beta|+1}\beta!^{\alpha} \quad \text{in } U,
$$
where $\beta!^{\alpha}= \beta_{0}!^{\alpha_{0}}\beta_{1}!^{\alpha_{1}}\dots \beta_{n}!^{\alpha_{n}}$.
\end{definition}
\medskip
Our result can be stated as follows:
\begin{proposition}\label{P-1}
Let $P_{1}$ be as in the Theorem \ref{T-1}, where $pk <\ell$.
If $u $ solves the problem 
$P_{1}u =f$ and $f$ is analytic then
$u \in G^{(s_{1},s_{2},s_{3},s_{4})}$ where
$s_{4} \geq \sup\lbrace\frac{r+kp}{q},\frac{r}{p}\rbrace$,
$s_{2} \geq \frac{k}{k+1} + \frac{1}{k+1} \frac{r+kp}{q}$,
$s_{3}\geq \frac{r(q-1)}{r(p-1)+q-p}$
and
$s_{1} \geq 1+ \sup \Big\lbrace \frac{1}{p(k+1)}\left( \frac{r+kp}{q} -1\right),  \frac{1}{r}\left( \frac{r+kp}{q} -1\right),
\frac{1}{r}\left(\frac{r}{p}-1\right), \frac{(r-1)(q-p)}{q(r(p-1)+q-p)}
\Big\rbrace$.
\end{proposition}
%
The same strategy used in the proof of the above Proposition shows that:
\begin{remark}
If $pk > \ell $ then
$u \in G^{(s_{1},s_{2},s_{3},s_{4})}$ where
$s_{4} \geq \sup\lbrace\frac{r+\ell}{q},\frac{r}{p}\rbrace$,
$s_{2} \geq \frac{k}{k+1} + \frac{1}{k+1} \frac{r+\ell}{q}$,
$s_{3}\geq \frac{r(q-1)}{r(p-1)+q-p}$
and
$s_{1} \geq 1 +\sup \Big\lbrace \frac{1}{p(k+1)}\left( \frac{r+\ell}{q} -1\right), \frac{1}{r}\left( \frac{r+\ell}{q} -1\right),
\frac{1}{r}\left(\frac{r}{p}-1\right), \frac{(r-1)(q-p)}{q(r(p-1)+q-p)}
\Big\rbrace$.
\end{remark}
%
\begin{remark}
Let $P_{2}(x;D)$ as in the Theorem \ref{T-1}.
If $qk < \ell $ and $u $ solves the problem $P_{2}u =f$, $f$ analytic, then
$u \in G^{(s_{1},s_{2},s_{3},s_{4})}$ where
$s_{4} \geq \frac{r+kq}{p}$, $s_{2}\geq 1$,
$s_{3} \geq \frac{k}{k+1} + \frac{1}{k+1} \frac{r+kq}{p}$
and
$s_{1} \geq \sup \Big\lbrace 1+\frac{1}{p(k+1)}\left( \frac{r+kq}{p} -1\right),  1+\frac{1}{r}\left( \frac{r+kq}{p} -1\right)
\Big\rbrace$. Otherwise
if $qk \geq \ell $ then 
$u \in G^{(s_{1},s_{2},s_{3},s_{4})}$ where
$s_{4} \geq \frac{r+\ell}{p}$, $s_{2}\geq 1$,
$s_{3} \geq \frac{k}{k+1} + \frac{1}{k+1} \frac{r+\ell}{p}$
and
$s_{1} \geq \sup \Big\lbrace 1+\frac{1}{p(k+1)}\left( \frac{r+\ell}{p} -1\right),  1+\frac{1}{r}\left( \frac{r+\ell}{p} -1\right)
\Big\rbrace$.
\end{remark}
\vspace{0.2cm}
\begin{proof}[\textbf{Proof Proposition} \text{\textbf{\ref{P-1}}}]
Since the regularity in the direction $D_{4}$ has been obtained in the Theorem \ref{T-1} we have only to analyze the 
direction $D_{1}$, $D_{2}$ and $D_{3}$. 
The primary tool will be once again the subelliptic estimate (\ref{SB-E-1/2}). Roughly speaking the strategy will be to transform the derivatives in the directions $D_{2}$ and $D_{1}$ in powers of the derivative in the direction $D_{4}$,
this will allow us to use the result in the Theorem \ref{T-1}.  Concerning the direction $D_{3}$ we will obtain the result directly.\\

{\bf Direction $ {\bm D_{3}}$}: Let $ \varphi_{N} (x_{3},x_{4})$ be a cut off function of Ehrenpreis-H\"or\-man\-der
type described in the proof of the Theorem \ref{T-1}-\textbf{\textit{A}} to analyze the direction $ x_{4} $.
We replay $ u $ by $ \varphi_{N} D_{3}^{N} u $ in ( \ref{SB-E-1/2}).
We have
\begin{equation}\label{RS-Est3}
\| \varphi_{N} D^{N}_{3} u\|^{2}_{\frac{1}{r+kp}} + \sum_{j=1}^{5} \|X_{j} \varphi_{N} D^{N}_{3} u \|^{2}_{0} 
\leq 
C \left( |\langle P_{1} \varphi_{N} D^{N}_{3} u, \varphi_{N} D^{N}_{3} u \rangle | + \| \varphi_{N} D^{N}_{3} u \|^{2}_{0} 
\right). 
\end{equation}
The scalar product in the right hand side leads to
\begin{align}\label{es-3.1} 
& 
2 \sum_{j=1}^{5} \langle X_{j} [X_{j}, \varphi_{N} ] D_{3}^{N} u,
 \varphi_{N} D_{3}^{N} u \rangle  
+   \sum_{j=1}^{5} \langle[ X_{j},[X_{j}, \varphi_{N} ] ] D_{3}^{N} u , 
 \varphi_{N} D_{3}^{N} u \rangle
\\
&\nonumber \qquad \qquad \qquad \qquad \qquad\qquad\qquad\qquad\qquad\qquad 
 + \langle \varphi_{N} D_{3}^{N} P_{1} u, \varphi_{N} D_{3}^{N} u \rangle.
\end{align}
The last term has a trivial estimate since $ P_{N}u $ is analytic. Without loss of generality
we can assume that it is zero. We focus our attention only on
the vector field $ X_{3} $, the case $ X_{4} $ and $ X_{5} $ can be handled
in the same way, these vector fields have coefficients with power of $ x_{1}$
greater than $ q-1 $. We have
\begin{align*}
&|\langle [X_{3}, \varphi_{N} ] D_{3}^{N} u, X_{3}\varphi_{N} D_{3}^{N} u \rangle |
=|\langle x_{1}^{q-1}\varphi^{(1)}_{N} D_{3}^{N} u, X_{3}\varphi_{N} D_{3}^{N} u \rangle |\\
&\qquad
\leq
|\langle X_{3} \varphi^{(1)}_{N}D_{3}^{N-1} u, X_{3}\varphi_{N} D_{3}^{N} u \rangle |
+
|\langle x_{1}^{q-1}\varphi^{(2)}_{N} D_{3}^{N-1} u, X_{3}\varphi_{N} D_{3}^{N} u \rangle |\\
&\qquad
\leq
C_{1} \| X_{3} \varphi^{(1)}_{N} D_{3}^{N-1} u\|^{2} 
+
\frac{1}{C_{1}} \| X_{3}\varphi_{N} D_{3}^{N} u \|^{2}
+
|\langle x_{1}^{q-1}\varphi^{(2)}_{N} D_{3}^{N-1} u, X_{3}\varphi_{N} D_{3}^{N} u \rangle |\\
&\qquad
\leq \dots \leq\quad
\sum_{j=1}^{N}C_{j}\| X_{3} \varphi^{(j)}_{N} D_{3}^{N-j} u\|^{2}
+
\frac{1}{C_{j}}\| X_{3}\varphi_{N} D_{3}^{N} u\|^{2}
+
C_{N+1}\|\varphi^{(N+1)}_{N} u \|^{2}\\
&\qquad\qquad\qquad\qquad\qquad\qquad\qquad\qquad\qquad\qquad\qquad\quad\quad
+
\frac{1}{C_{N+1}} \| X_{3}\varphi_{N} D_{3}^{N} u\|^{2}.
\end{align*}
The constant $ C_{j}$ are arbitrary, we make the choice $ C_{j} = \delta^{-1} 2^{j}$,
for a suitable fixed small $ \delta $. We can absorb each term of the form 
$ C^{-1}_{j} \| X_{3} \varphi_{N} D^{\alpha}_{3} u\|^{2}$ on the left hand side of
(\ref{RS-Est3}).
The term $ C_{N+1} \| \varphi^{N+1}_{N} u \|^{2}$ is smaller than $ C^{2(N+1) } N!^{2} $,
that is it gives analytic growth.
To estimate the terms $ C_{j}\|X_{3} \varphi^{j}_{N} D_{3}^{N-j} u\|^{2}$, we observe
that for each of them there has been a shift of one or more $x_{3}$-derivatives
from $ u $ to $ \varphi_{N} $, but they have the same form as $\| X_{3} \varphi_{N} D_{3}^{N} u \|^{2}$.
We have to estimate the sum
\begin{align}\label{sum-3}
\sum_{j=1}^{N} \frac{2^{j}}{\delta} &\| X_{3}\varphi^{(j)}_{N} D_{3}^{N-j}u\|^{2}
=
\frac{2}{\delta} \| X_{3}\varphi^{(1)}_{N} D_{3}^{N-1}u\|^{2}
+
\sum_{j=2}^{N} \frac{2^{j}}{\delta} \| X_{3}\varphi^{(j)}_{N} D_{3}^{N-j}u\|^{2}.
\end{align}
We start from the first term in the sum. We use the Rothschild-Stein
sub-elliptic estimate replacing $ u$ with $ \varphi^{(1)}_{N} D_{3}^{N-1} u$,
repeating the above procedure we have
\begin{align*}
\| X_{3} \varphi^{(1)}_{N}D_{3}^{N-1} u \|^{2}
\leq 
\sum_{j=1}^{N-1} 
\left( \frac{2^{j}}{\delta}\| X_{3} \varphi^{(j+1)}_{N} D_{3}^{N-j-1} u\|^{2} 
+ \frac{\delta}{2^{j}} \| X_{3} \varphi^{(1)}_{N}D_{3}^{N -1} u\|^{2} 
\right) 
\end{align*}
modulo terms  which give analytic growth or which have the following form 
$|\langle  [ X_{3},[ X_{3}, \varphi^{(1)}_{N} ] ] D_{3}^{N -1} u, \varphi^{(1)}_{N} D_{3}^{N-1} u \rangle | $;
we remark that for each of them there has been a shift of $ x_{3}$-derivatives
from $ u $ to $ \varphi_{N} $, but essentially they have the same form as
$|\langle  [ X_{3},[ X_{3}, \varphi_{N} ] ] D_{3}^{N } u, \varphi^{(1)}_{N} D_{3}^{N-1} u \rangle | $
in (\ref{es-3.1}), for the discussion of these terms see in the continuations of the proof.
As before we may absorb the second term in the left hand side of the estimate.
Repeating the above process $s$ times we have
\begin{equation*}
\sum_{j =1}^{N} \frac{1}{\delta}2^{j}
\| X_{3}\varphi^{(j)}_{N} D_{3}^{N -j} u\|^{2}
\leq 
\sum_{j =s}^{N} \frac{1}{\delta}\left(1+\frac{1}{\delta}\right)^{s-1} 2^{j}
\| X_{3}\varphi^{(j)}_{N} D_{3}^{N -j} u\|^{2}
\end{equation*}
modulo terms which can be absorbed on the left hand side or which
give analytic growth or which have the form
$|\langle  [ X_{3},[ X_{3}, \varphi^{(j)}_{N} ] ] D_{3}^{N -1} u, \varphi^{(j)}_{N} D_{3}^{N-j} u \rangle | $,
$ 1\leq j \leq s-1$.
With the same procedure, after $ N -1$ iterates, we obtain a term of the form
\begin{equation*}
\frac{1}{\delta}\left( 1+\frac{1}{\delta}\right)^{N-1}
2^{N} \| X_{3} \varphi^{(N)}_{N} u\|^{2}.
\end{equation*}
This term can be estimate
by $ C^{2(N+1)} (N!)^{2}$, we have analytic growth.\\
\noindent
On the other hand we have
\begin{align*}
|\langle  [ X_{3},[ X_{3}, \varphi_{N} ] ] D_{3}^{N } u, \varphi^{(1)}_{N} D_{3}^{N-1} u \rangle | 
&=
| \langle x_{1}^{2(q-1)} \varphi^{(2)}_{N}D_{3}^{N} u,\varphi_{N} D_{3}^{N} u \rangle |\\
&\leq
| \langle x_{1}^{q-1} D_{3}\varphi^{(2)}_{N}D_{3}^{N-1} u, x_{1}^{q-1} D_{3}\varphi_{N} D_{3}^{N-1} u \rangle |\\
&\quad
+| \langle x_{1}^{q-1} D_{3}\varphi^{(2)}_{N}D_{3}^{N-1} u, x_{1}^{q-1} \varphi^{(1)}_{N} D_{3}^{N-1} u \rangle |\\
&\quad
+| \langle x_{1}^{q-1} \varphi^{(3)}_{N}D_{3}^{N-1} u, x_{1}^{q-1}D_{3} \varphi_{N} D_{3}^{N-1} u \rangle |\\
&\quad
+| \langle x_{1}^{q-1}\varphi^{(3)}_{N}D_{3}^{N-1} u, x_{1}^{q-1} \varphi^{(1)}_{N} D_{3}^{N-1} u \rangle |\\
&
= H_{0}+H_{1}+H_{2}+H_{3}.
\end{align*}
We study any single term. Term $ H_{0}$:
\begin{align*}
H_{0} \leq \frac{2}{N^{2}} \| X_{3} \varphi^{(2)}_{N} D_{3}^{N-1} u \|^{2}
+ 2 N^{2} \| X_{3} \varphi_{N} D_{3}^{N-1} u \|^{2}.
\end{align*}
As done previously the weight $ N $ is introduced to balance the number of
$ x_{3}$-derivatives on $ u$ with the number of derivatives on $ \varphi_{N}$.
The terms on the right hand side have the same form as $ \| X_{3} \varphi_{N} D_{3}^{N} u \|^{2}$.
We can restart the process.\\
\noindent
The term $ H_{1}$:
\begin{align*}
H_{1} 
&\leq C_{1} \| X_{3} \varphi^{(2)}_{N} D_{3}^{N -2} u \|^{2} 
+ \frac{1}{C_{1}} \| X_{3} \varphi^{(1)}_{N} D_{3}^{N -1} u \|^{2}\\
&\qquad\qquad\qquad\qquad\qquad
+ | \langle x_{1}^{q-1}\varphi^{(3)}_{N}D_{3}^{N-2} u, x_{1}^{q-1}D_{3} \varphi^{(1)}_{N} D_{3}^{N-1} u \rangle |\\
&\leq \dots \leq
\sum_{j=1}^{N} \left( C_{j} \| X_{3} \varphi^{(j)}_{N} D_{3}^{N-j} u \|^{2}
+ \frac{1}{C_{j}}  \| X_{3} \varphi^{(1)}_{N} D_{3}^{N-1} u \|^{2}\right)
+ C_{N} \| \varphi^{(N+1)}_{N} u \|^{2}.
\end{align*}
The above sum can be handled with the same strategy used to estimate the sum (\ref{sum-3}).
The last term give analytic growth.\\
\noindent
The term $ H_{2}$:
\begin{align*}
H_{2} &\leq 
\frac{C_{1}}{N^{4}}\| X_{3} \varphi^{(3)}_{N}  D_{3}^{N-1}  u \|^{2} + 
\frac{N^{4}}{C_{1}}\| X_{3} \varphi_{N}  D_{3}^{N-2}  u \|^{2}+
\frac{C_{2}}{N^{2}} \| X_{3} \varphi^{(3)}_{N}  D_{3}^{N-2}  u \|^{2}  
\\
&\quad 
+\frac{N^{2}}{C_{2}} \| X_{3} \varphi^{(1)}_{N}  D_{3}^{N-2}  u \|^{2}+
|\langle x_{1}^{q-1} \varphi^{(4)}_{N} D_{3}^{N-2} u, x_{1}^{q-1}D_{3} \varphi^{(1)}_{N} D_{3}^{N-2}u \rangle|
\\
& 
\leq \ldots \leq  
\sum_{j=1}^{N} \frac{C_{1}}{N^{4}}\| X_{3}\varphi^{(j+2)}_{N}D_{3}^{N-j}u\|^{2} 
+\sum_{j=1}^{N -1}\frac{N^{4}}{C_{1}} \| X_{3}\varphi^{(j-1)}_{N}D_{3}^{N-(j+1)}u\|^{2}\\ 
&\qquad\qquad
+\sum_{j=1}^{N} \left\{ \frac{C_{2}}{N^{2}}\| X_{3}\varphi^{(j+2)}_{N}D_{3}^{N-(j+1)}u\|^{2}
+\frac{N^{2}}{C_{2}}\| X_{3}\varphi^{(j)}_{N}D_{3}^{N-(j+1)}u\|^{2}\right\}\\ 
& \qquad\qquad
+ | \langle x_{1}^{q-1} \varphi^{(N+2)}  u, x_{1}^{q-1} D_{3}\varphi^{(N -1)}_{N} u \rangle|.
\end{align*}
%
The last term gives analytic growth.
To estimate the terms in the sums, we observe that  with the help of the weight $N$ we have essentially, 
on each of them,  shifted  one or more $ x_{3} $-derivatives from $ u $ to $ \varphi_{N} $;
they have the same form as $ \| X_{3} \varphi_{N} D_{3}^{N} u\|^{2} $.\\
\noindent
The term $H_{3}$:
\begin{align*}
H_{3} & \leq
|\langle x_{1}^{q-1}D_{3} \varphi^{(3)}_{N} D_{3}^{N-2} u, x_{1}^{q-1}D_{3} \varphi^{(1)}_{N} D_{3}^{N-2}u \rangle|\\ 
& \quad 
+
|\langle x_{1}^{q-1}D_{3} \varphi^{(3)}_{N} D_{3}^{N-2} u, x_{1}^{q-1}D_{3} \varphi^{(1)}_{N} D_{3}^{N-2}u \rangle|\\
&\quad 
+
|\langle x_{1}^{q-1} \varphi^{(4)}_{N} D_{3}^{N-2} u, x_{1}^{q-1} D_{3}\varphi^{(1)}_{N} D_{3}^{N-2} u \rangle|\\
&\quad 
+
 |\langle x_{1}^{q-1} \varphi^{(4)}_{N} D_{3}^{N-2} u, x_{1}^{q-1} \varphi^{(2)}_{N} D_{3}^{N-2} u \rangle|. 
\end{align*}
Iterating we obtain
\begin{align*}
H_{3} 
&\leq 
\sum_{j=1} 
|\langle x_{1}^{q-1}D_{3} \varphi^{(j+2)}_{N} D_{3}^{N-(j+1)} u, x_{1}^{q-1}D_{3} \varphi^{(j)}_{N} D_{3}^{N-(j+1)} u \rangle|\\ 
&\quad
+
\sum_{j=1}  
| \langle x_{1}^{q-1} \varphi^{(j+2)}_{N} D_{3}^{N-(j+1)} u, x_{1}^{q-1} \varphi^{(j+1)}_{N} D_{3}^{N-(j+1)} u \rangle|\\
&\quad
+
\sum_{j=1}
| \langle x_{1}^{q-1} \varphi^{(j+3)}_{N} D_{3}^{N-(j+1)} u, x_{1}^{q-1} D_{3}\varphi^{(j)}_{N} D_{3}^{N-(j+1)}u \rangle|\\
&\quad
+ 
| \langle x_{1}^{q-1} \varphi^{(N+2)}_{N} u, x_{1}^{q-1} \varphi^{(N)}_{N} u \rangle| . 
\end{align*}
We observe that the terms in the first sum have the same form as $H_{0}$, 
the terms in the second sum have the same form as $H_{1}$ and
those in the third  sum have the same form as $ H_{2}$, 
we can handle each of them as above.
Finally, the last term gives analytic growth.
Using the estimate (\ref{SB-E-1/2}) with $ u $ replaced by $N^{i}\varphi^{(j)}_{N} D_{3}^{N-(j+i)} u $ 
or  $ N^{-i}\varphi^{(j+i)}_{N} D_{3}^{N -j} u $ and applying recursively  the same strategy followed above
we are able to shift all free derivatives on $ \varphi_{N} $.\\
\noindent
As previously observed, to analyze the case $ X_{4}$ and $ X_{5}$ we can use the same 
strategy used to study the case $ X_{3}$. Indeed since the commutators
$ [ X_{4}, \varphi_{N}]$, $ [ X_{5}, \varphi_{N}]$, $[X_{4}, [ X_{4}, \varphi_{N}]]$
and $ [X_{5},[ X_{5}, \varphi_{N}]]$ give terms with powers of $ x_{1}$
greater than $ q-1$, we can take again maximum advantage
from the sub-elliptic estimate. Also in these cases we have 
analytic growth.\\ 
Hence we have
\begin{equation*}
\| \varphi D_{3}^{N} u \|^{2}_{\frac{1}{r+kp}} 
+  
\sum_{j=1}^{5} \| X_{j} \varphi_{N} D_{3}^{N} u \|^{2}
\leq 
C^{2(N+1)} N^{2N}.
\end{equation*}
To obtain the result we need to consider when $x_{2} \neq 0$.
To do it since when $x_{2}\neq 0$ the operator $P_{1}$ is an operator
of Ole\u{\i}nik-Radkevi\v c type, \cite{OR1973}, we use the following result in \cite{C}: 
\begin{theorem}[\cite{C}]\label{st-C}
Let $P$ be the operator given by
\begin{equation}\label{OperatorL}
P(x,D_{x}) =
D_{x_{1}}^{2} + \sum_{j=2}^{n} x_{1}^{2(r_{j}-1)}D_{x_j}^{2}.
\end{equation}
We have that $P$ is $G^{r_{n}/r_{1}}$ hypoelliptic and not better. More precisely 
we have that if $u$ solves the equation $Pu=f$ and $f$ is analytic then if $\rho_{j}\in \Sigma_{r_{j}-1}$
then $\rho_{j}\notin WF_{r_{j}/r_{1}}(u)$ and moreover $u \in G^{(s_{0},s_{1},\dots,s_{n})}$ where
%
$$
s_{1}\geq r^{*},\quad  s_{j} = \beta_{j} \geq \frac{r_{n}(r_{j}-1)}{r_{n}(r_{1}-1)+r_{j}-r_{1}}
\text{ with } j=2,\dots,n ;
$$
where $ r^{*} = \displaystyle\sup_{j}\left\{ 1-\frac{1}{r_{j}}+\frac{\beta_{j}}{r_{j}}\right\} $,
in particular $s_{2}\geq 1$ and $s_{n}\geq r_{n}/r_{1}$.
\end{theorem}
We can conclude that 
we have in the direction $x_{3}$ a growth corresponding
to $ G^{\frac{r(q-1)}{r(p-1)+q-p}}$. 


\vspace*{0.4cm}

{\bf Direction $ {\bm D_{2}}$.} Once again our primary tool will be the sub-elliptic estimate
(\ref{SB-E-1/2}). As in the study of the direction $ x_{3}$,
we replace $ u $ by $ \varphi_{N} D_{2}^{N} u $ in ( \ref{SB-E-1/2}).
We recall that $ \varphi_{N} $ does not depend on $ x_{1}$ and $ x_{2}$. We have
\begin{equation}
\label{RS-Est2}
\| \varphi_{N} D^{N}_{2} u\|^{2}_{\frac{1}{r+kp}} \!\!+\! \sum_{j=0}^{5} \|X_{j} \varphi_{N} D^{N}_{2} u \|^{2}_{0} 
\leq 
C\! \left( |\langle P_{1} \varphi_{N} D^{N}_{2} u, \varphi_{N} D^{N}_{2} u \rangle | + \| \varphi_{N} D^{N}_{2} u \|^{2}_{0} 
\right). 
\end{equation}
We consider the scalar product in the right hand side of the above inequality.
We have to study terms of the type
\begin{equation*}
| \langle [ X_{j}, \varphi_{N} D_{2}^{N}] u , X_{j} \varphi_{N} D_{2}^{N} u\rangle |
\quad j = 3,4, 5.
\end{equation*}
%
Since $ X_{3} =x^{q-1}_{1} D_{3} $, $ X_{5} = x^{r +\ell -1}_{1} D_{4}$
and $ q $  and $ r $ are strictly greater than $p $, as seen in the study
of the direction $ x_{3}$, we can take maximum advantage from the sub-elliptic
estimate shifting one derivative from $ u$ to $ \varphi_{N}$. If we focus our attention
only on these terms and we iterate the process we will obtain analytic growth.\\
\noindent
The case $ X_{4}= x_{1}^{r-1} x_{2}^{k} D_{4}$. We have
\begin{align*}
[ X_{4}, \varphi_{N} &D_{2}^{N}] u
= 
[ x_{1}^{r-1}x_{2}^{k} D_{4},  \varphi_{N} D_{2}^{N}] u
=   x_{1}^{r-1} x^{k}_{2}\varphi^{(1)}_{N} D_{2}^{N} u 
+ x_{1}^{r-1} \varphi_{N} [x_{2}^{k}, D_{2}^{N}] D_{4} u\\
&
= x_{1}^{r-1} x^{k}_{2}\varphi^{(1)}_{N} D_{2}^{N} u
- x_{1}^{r-1} \varphi_{N} \sum_{j = 1}^{k} 
\frac{N! k!}{ (i)^{j} j! ( N -j)! (k-j)!} x_{2}^{k-j} D_{2}^{N -j} D_{4} u.
\end{align*}
Without loss of generality we analyze one of the terms; a similar method 
can be used to handle the other terms. We consider the first one:
$ N k i^{-1}  x_{1}^{r-1}D_{2} \varphi_{N} x_{2}^{k-1} D_{2}^{N -2} D_{4} u$.
We have to estimate 
$ N k  \| X_{2} \varphi_{N} x_{2}^{k-1} D_{2}^{N -2} D_{4} u \|$.
We apply the sub-elliptic estimate with 
$ u $ replayed by $N k  \varphi_{N} x_{2}^{k-1} D_{2}^{N -2} D_{4} u$, arguing as above,
we study the first term coming from the commutator with $ X_{4}$.
We obtain the term 
$ k^{2} N (N -2) x_{1}^{r-1} x_{2}^{2(k-1)} \varphi_{N} D_{2}^{N -3} D_{4}^{2}u$.
We have to estimate 
$ k^{2} N (N -2) \| X_{4} x_{2}^{k-2} \varphi_{N} D_{2}^{N -3} D_{4} u \| $.
Hence after two steps we have
\begin{equation*}
\| X_{4} \varphi_{N} D_{2}^{N} u \| 
\rightarrow
k^{2} N (N -2) \| X_{4} x_{2}^{k-2} \varphi_{N} D_{2}^{N -3} D_{4}u \|. 
\end{equation*}
Repeating the process $ j$-times, we have
\begin{equation*}
\| X_{4} \varphi_{N} D_{2}^{N} u \| 
\rightarrow \dots  \rightarrow
C \frac{N!}{(N -1)( N- (j+1))!}  \| X_{4} x_{2}^{k-j} \varphi_{N} D_{2}^{N -(j+1)} D_{4}u \|. 
\end{equation*}
Here the constant $ C $ depend by $ k$.
We stress that $ N![(N -1)( N- (j+1))!]^{-1} \sim N^{j} $.
In this way after $ k $ iterates we have to analyze a term of the form
$ C_{k}  N![(N -1)( N- (k+1))!]^{-1} \varphi_{N} D_{2}^{N - (k+1)} D_{4} u$.
Arguing in the same way after $ m $ steps we have
\begin{equation*}
\| X_{4} \varphi_{N} D_{2}^{N} u \| 
\rightarrow \dots  \rightarrow
C_{k} 
N^{m k }
\| X_{4}  \varphi_{N} D_{2}^{N -m(k+1)} D_{4}^{m}u \|. 
\end{equation*}
Iterating the cycle $N /(k+1)$-times we use up all free derivatives in
$ x_{2}$-direction and we are left with
\begin{equation*}
C_{k}^{N}   
N^{N\frac{k}{k+1}}
\| X_{4}  \varphi_{N}  D_{4}^{\frac{N}{k+1}}u \|. 
\end{equation*}
As well as it was done in the proof of the Theorem \ref{T-1} we introduce
$\chi_{N}(\xi_{4})$ an Ehrenpreis-H\"ormander cutoff function such that
$\chi_{N}$ is $ C^{\infty}(\mathbb{R})$ non negative function such that $\chi_{N} = 0$
for $\xi_{4} < 3 $ and $\chi_{N} =1$ for $\xi_{4} > 4$.
We have
\begin{equation*}
\| X_{4}  \varphi_{N} \chi_{N}(N^{-1}\xi_{4}) D_{4}^{\frac{N}{k+1}}u \| 
\leq \| X_{4} \varphi_{N} \left( 1-\chi_{N}(N^{-1}D_{4})\right)D_{4}^{\frac{N}{k+1}} u \|
+ \| X_{4} \varphi_{N} \chi_{N}(N^{-1}D_{4})D_{4}^{\frac{N}{k+1}} u \|.
\end{equation*} 
Since $ 1-\chi_{N}(N^{-1}D_{4})$ has support for $\xi_{4} < 4N$ we have
\begin{equation*}
C_{k}^{N}   
N^{N\frac{k}{k+1}}\| X_{4}\varphi_{N} \left( 1-\chi_{N}(N^{-1}D_{4})\right)D_{4}^{\frac{N}{k+1}} u \|
\leq C^{N+1}N^{N},
\end{equation*} 
where $C$ is a positive constant independent by $N$, but depending on $u$ and $k$.
To estimate $\| X_{4} \varphi_{N} \chi_{N}(N^{-1}D_{4})D_{4}^{\frac{N}{k+1}} u \|$ we use the
same strategy used in the proof of the Theorem \ref{T-1}.
Therefore since in the direction $ x_{4}$ we have a growth corresponding to
$ G^{\frac{r +kp}{q}}$ we can estimate this term with
$ C^{N +1} (N!)^{\frac{r+k(p+q)}{q(k+1)}}$. 
We can estimate the left hand side of (\ref{RS-Est2}) with this quantity, we have the growth
corresponding to $ G^{(r+k(p+q))/q(k+1)} $.\\
\noindent
More in general 
applying the sub-elliptic estimate  
and iterating the above processes more time, we may estimate 
the left hand side of (\ref{RS-Est2}) with terms of the form
\begin{equation*}
(N)^{(N-j)mk} \| X_{4} \varphi^{(j)}_{N} D_{2}^{N -j -m(k+1)}D_{4}^{N} u \|.
\end{equation*}
Iterating the procedure until all the $ x_{2}$-derivatives are used up we have 
to apply the sub-elliptic estimate to terms of the form
\begin{equation*}
(N)^{(N-j)\frac{k}{k+1}}  \varphi^{(j)} D_{4}^{\frac{N -j}{k+1}}u .
\end{equation*}
To handle these terms we argue as before that is we introduce the cut-off $\chi_{N}$ 
and we apply the strategy used in the proof of the Theorem \ref{T-1} to obtain the
Gevrey regularity in the direction $ x_{4}$. 
Since  $ (r+kp)/q  >1$ we can conclude
\begin{equation*}
\| \varphi_{N} D^{N}_{2} u\|^{2}_{\frac{1}{r+kp}} + \sum_{j=0}^{5} \|X_{j} \varphi_{N} D^{N}_{2} u \|^{2}_{0} 
\leq 
C^{N +1} 
(N!)^{\frac{1}{k+1}\left(\frac{r+kp}{q}+k \right)}.
\end{equation*}
To gain the result we  need to consider when $x_{2} \neq 0$.
To do it since when $x_{2}\neq 0$ the operator $P_{1}$ is an operator
of Ole\u{\i}nik-Radkevi\v c type, \cite{OR1973}, we use Theorem \ref{st-C}.
We have that when $x_{2}\neq 0$ in the direction $D_{2}$ we have analytic growth.
We conclude that in this direction the growth corresponding to $G^{(r+k(p+q))/q(k+1)}$.\\
\vspace*{0.4cm}

{\bf Direction $ {\bm D_{1}}$}:
As in the study of the other directions,
we replace $ u $ by $ \varphi_{N} (x) D_{1}^{N} u $ in ( \ref{SB-E-1/2}).
We have
\begin{equation}
\label{RS-Est1}
\| \varphi_{N} D^{N}_{1} u\|^{2}_{\frac{1}{r+kp}} \!\!+ \!\sum_{j=0}^{5} \|X_{j} \varphi_{N} D^{N}_{1} u \|^{2}_{0} 
\leq 
C\! \left( |\langle P_{1} \varphi_{N} D^{N}_{1} u, \varphi_{N} D^{N}_{1} u \rangle | + \| \varphi_{N} D^{N}_{1} u \|^{2}_{0} 
\right). 
\end{equation}
We consider the scalar product in the right hand side of the above inequality.
We have to study terms of the type
\begin{equation*}
| \langle [ X_{j}, \varphi_{N} D_{1}^{N}] u , X_{j} \varphi_{N} D_{1}^{N} u\rangle |, \quad j = 2, 3, 4, 5.
\end{equation*}
%
We describe the case $ X_{2}$, the other cases can be handled using the same strategy.
We have
\begin{align*}
[ X_{2}, \varphi_{N} D_{1}^{N}] u
=  \varphi_{N} \sum_{j = 1}^{p-1} 
\frac{N! (p-1)!}{ (i)^{j} j! ( N -j)! (p-1-j)!} x_{1}^{p-1-j} D_{1}^{N -j} D_{2} u.
\end{align*}
Without loss of generality we analyze one of the terms. A similar method can be used
to handle the other terms.
Consider  $ N (p-1) D_{1} x_{1}^{p-2}\varphi_{N}D_{1}^{N -2} D_{2} u $ that is we have to estimate
a term of the form $N (p-1) \| X_{1} x_{1}^{p-2}\varphi_{N}D_{1}^{N -2} D_{2} u \|$. 
Applying the sub-elliptic estimate with $ u$ replaced by $ x_{1}^{p-2}\varphi_{N}D_{1}^{N -2} D_{2} u$
and arguing as above, we study the first term coming from the commutator
with $ X_{2}$. We obtain the term $ N (N -2)(p-1)^{2}  x_{1}^{2(p-2)}D_{1}^{N -3} D_{2}^{2} u $.
We have to estimate $ N (N -2)(p-1)^{2}\|X_{2}  x_{1}^{p-3}D_{1}^{N -3} D_{2} u\| $.
Hence after two step we have
\begin{equation*}
\| X_{2} \varphi_{N} D_{1}^{N} u \| 
\rightarrow
(p-1)^{2} \frac{N!}{(N -1)(N-3)!} \| X_{2} x_{1}^{p-3} \varphi_{N} D_{1}^{N -3} D_{2}u \|. 
\end{equation*}
Repeating the process $ s$-times, we have
\begin{equation*}
\| X_{2} \varphi_{N} D_{2}^{N} u \| 
\rightarrow \dots  \rightarrow
C_{p} \frac{N!}{(N -1)( N- (s+1))!}  \| X_{2} x_{1}^{p-(s+1)} \varphi_{N} D_{1}^{N -(s+1)} D_{2}u \|. 
\end{equation*}
We stress that $ N![(N -1)( N- (j+1))!]^{-1} \sim N^{j} $.
In this way after $ s = p-1$ iterates we have to analyze a term of the form
$ C_{p}  N^{p-1}\| X_{2} \varphi_{N} D_{1}^{N - p} D_{2} u \|$.
Arguing in the same way after $ m $ steps we have   
\begin{equation*}
\| X_{2} \varphi_{N} D_{1}^{N} u \| 
\rightarrow \dots  \rightarrow
C_{p}^{m} 
N^{m(p-1)}
\| X_{2}  \varphi_{N} D_{1}^{N -mp} D_{2}^{m}u \|. 
\end{equation*}
Iterating the cycle $N /p$-times we use up all free derivatives in
$ x_{1}$-direction and we are left with
\begin{equation*}
C_{p}^{N}   
N^{N\left( 1 - \frac{1}{p}\right)}
\| X_{2}  \varphi_{N}  D_{2}^{\frac{N}{p}}u \|. 
\end{equation*}
Since in the direction $ x_{2}$ we have a growth as 
$ G^{\frac{r +k(p+q)}{q(k+1)}}$ we can estimate the above term with
$$
C^{N +1} (N!)^{1 + \frac{1}{p}\left( \frac{r+k(p+q)}{q(k+1)}-\frac{1}{k+1}\right)}.
$$
We have the growth
$ G^{1 +  \frac{r+kp-q}{pq(k+1)}} $.\\
The other cases, that is the terms involving the commutators
with $ X_{3} $, $ X_{4}$ and $ X_{5}$, can be handled in the 
same way achieving  analytic growth, 
$ 1 + (r + kp-q)/rq$-Gevrey growth and  $ 1 + (r + kp-q)/(r+\ell)q$-Gevrey growth
respectively. We remark that in these three situations,
arguing as above, we obtain terms of the form
$ C_{q}^{N} (N ! )^{(q-1)/q} \| X_{3} \varphi_{N} D_{3}^{N/q} u\| $,
$ C_{r}^{N} (N ! )^{(r-1)/r} \| X_{4} \varphi_{N} D_{4}^{N/r} u\| $
and
$ C_{r+\ell}^{N} (N ! )^{(r + \ell -1)/( r + \ell )} \| X_{5} \varphi_{N} D_{4}^{N/( r + \ell)} u\| $.
Moreover we point out that also in the general situation
we will obtain a Gevrey growth  less than or equal to that obtained by analyzing the individual cases.
We have obtained a growth corresponding to $ G^{s_{1}} $ 
where $ s_{1} =  \sup \lbrace 1 + \frac{r+kp -q}{qp(k+1)} , 1 + \frac{r + kp-q}{rq}\rbrace$.\\
To obtain the result we  need to consider when $x_{2} \neq 0$.
To do it since when $x_{2}\neq 0$ the operator $P_{1}$ is an operator
of Ole\u{\i}nik-Radkevi\v c type, \cite{OR1973}, we use Theorem \ref{st-C}.
We have that when $x_{2}\neq 0$ in the direction $D_{2}$ we have a 
growth corresponding to $ G^{s_{2}} $ where
$ s_{2} =  \sup \lbrace 1 + \frac{1}{q}\left(\frac{r(q-1)}{r(p-1)+q-p} -1\right),
1 + \frac{1}{r}\left(\frac{1}{p}-\frac{r}{p}\right)\rbrace$.
We conclude that in the direction $x_{2}$ we have a growth corresponding to
$G^{s}$ where $s =\sup \lbrace s_{1},s_{2}\rbrace$.
We point out that the case $x_{2} \neq 0$ can be directly considered taking
the cutoff function $\varphi_{N}$ depending also on the $x_{2}$-variable from the beginning.\\
\end{proof}
%
%
\section{Additional material: the $\textit{\textbf{n}}-$dimensional case}
Following the some ideas used to archive the Theorems \ref{T-1} and \ref{T-2} we can extend
without particular difficulties such results to the following $n$-dimensional cases, $n \geq 5$. We omit the proofs. 
%
\begin{theorem}\label{T-3}
Let $ P_{i,n}(x;D)$ be the operator given by
\begin{equation}\label{Op_n_i}
P_{i,n}(x;D) = 
D_{1}^{2} +\sum_{j=2}^{n-1} x_{1}^{2(r_{j}-1)} D_{j}^{2} 
+ \left( x_{1}^{2(r_{n}-1)} x_{i}^{2k} + x_{1}^{2(r_{n}+ \ell -1 )} \right)D_{n}^{2},
\qquad 2\leq i\leq n-1,
\end{equation}
in $\Omega$, open neighborhood of the origin in $\mathbb{R}^{n}$,
where $ r_{j}$, $j= 1,\dots, n$, $k$ and $\ell$ are positive integers such that $r_{1}< r_{2} < \dots < r_{n}$.
We have:
\vspace{-0.3cm}
\begin{itemize}
\item[i)] if $kr_{i} <\ell$, $P_{i,n}(x;D)$ is $G^{s}$-hypoelliptic with 
             $s= \displaystyle\sup\Big\lbrace\frac{r_{n}+kr_{2}}{r_{3}}, \frac{r_{n}}{r_{2}}\Big\rbrace$ if $i =2$
             and $s= \displaystyle\frac{r_{n}+kr_{i}}{r_{2}}$ if $i \neq 2$. In particular if $u$ solves the equation
             $P_{i,n}u = f$ and $f$ is analytic then the point
             $(0,e_{n})$ in $\Char(P_{i,n})$ does not belong to $\WF_{(r_{n}+kr_{2})/r_{3}}(u)$ if $i=2$ and
             it  does not belong to $\WF_{(r_{n}+kr_{i})/r_{2}}(u)$ if $i\neq2$.
             \vspace{-0.3cm}
\item[ii)]if $kr_{i} \geq\ell$, $P_{i,n}(x;D)$ is $G^{s}$-hypoelliptic with 
             $s= \displaystyle\sup\Big\lbrace\frac{r_{n}+\ell}{r_{3}}, \frac{r_{n}}{r_{2}}\Big\rbrace$ if $i =2$
             and $s= \displaystyle\frac{r_{n}+\ell}{r_{2}}$ if $i \neq 2$. In particular if $u$ solves the equation
             $P_{i,n}u = f$ and $f$ is analytic then the point
             $(0,e_{n})$ in $\Char(P_{i,n})$ does not belong to $\WF_{(r_{n}+\ell)/r_{3}}(u)$ if $i=2$ and
             it  does not belong to $\WF_{(r_{n}+\ell)/r_{2}}(u)$ if $i\neq2$.
\end{itemize}
\end{theorem}
%
\begin{remark}
Let $ \widetilde{P}_{i,n}(x;D)$ be the operator given by
\begin{equation}\label{Op_n_i-v}
\widetilde{P}_{i,n}(x;D) = 
D_{1}^{2} +\sum_{j=2}^{m-1} x_{1}^{2(r_{j}-1)} D_{j}^{2} 
+ \left( x_{1}^{2(r_{m}-1)} x_{i}^{2k} + x_{1}^{2(r_{m}+ \ell -1 )} \right)D_{m}^{2}
+\sum_{j=m+1}^{n} x_{1}^{2(r_{j}-1)} D_{j}^{2},
\end{equation}
in $\Omega$, open neighborhood of the origin in $\mathbb{R}^{n}$,
where $m\geq 3$, $2\leq i\leq m-1$, $ r_{j}$, $j= 1,\dots, n$, $k$ and $\ell$ positive integers such that 
$r_{1}< r_{2} < \dots < r_{n}$ and $r_{n}> r_{m}+\sup\lbrace kr_{i}, \ell\rbrace$. We have that $\widetilde{P}_{i,n}(x;D)$
is $r_{n}/r_{2}$-Gevrey hypoelliptic. In particular if $u$ solves the equation
$\widetilde{P}_{2,n}(x;D)u = f$ and $f$ is analytic then
the point $(0,e_{n})\in \Char(\widetilde{P}_{2,n})$ does not belong to $\WF_{r_{n}/r_{3}}(u)$.
\end{remark}
%
%
\begin{theorem}\label{T-4}
Let $P_{m,n}(x,D)$ be the operator given by
\begin{equation}\label{Op-m_n}
P_{m,n}(x;D) = D_{1}^{2} +\! \sum_{i=2}^{m+1} x_{1}^{2(r_{i}-1)} D_{i}^{2} + \!\!
\sum_{i=m+2}^{n}\!\!\left( x_{1}^{2(r_{i}-1)}x_{i-m}^{2k_{i-m}} \!\!+ x_{1}^{2(r_{i}+\ell_{i-m}-1)} \right) D_{i}^{2},
\,\,\, m \geq  \left[\frac{n}{2}\right],
\end{equation}
in $\Omega $, open neighborhood of the origin in $\mathbb{R}^{n}$, where
$r_{i}$, $i=2,\dots, n$, $k_{i-m}$ and $\ell_{i-m}$, $i=m+2,\dots, n$, are positive
integers such that $r_{2}< \cdots<r_{n}$, $k_{2}< \cdots<k_{n-m}$, $\ell_{2}< \cdots<\ell_{n-m}$ and
$r_{i}k_{i-m} < \ell_{i-m}$ for every $i$, $i = m+2,\dots, n$, then
$P_{m,n}(x,D)$ is $G^{s}$-hypoelliptic with $s= \displaystyle\frac{r_{n}+r_{n-m}k_{n-m}}{r_{2}}$.
Moreover if $u$ solves the equation $P_{m,n}u = f$ and $f$ is analytic then 
the point $(0,e_{n})\in \Char(P_{m,n})$ does not belong to  $WF_{\frac{r_{n}+r_{n-m}k_{n-m}}{r_{n-m+1}}}(u)$.
\end{theorem}
%
%
\begin{remark}
Let $\widetilde{P}_{m,n}(x,D)$ be the operator given by
\begin{equation*}
\widetilde{P}_{m,n}(x;D) = D_{1}^{2} +\! \sum_{i=2}^{m+2} x_{1}^{2(r_{i}-1)} D_{i}^{2} + \!\!
\sum_{i=m+3}^{n}\!\!\left( x_{1}^{2(r_{i}-1)}x_{i-m}^{2k_{i-m}} \!\!+ x_{1}^{2(r_{i}+\ell_{i-m}-1)} \right) D_{i}^{2},
\,\,\, m \geq  \left[\frac{n}{2}\right],
\end{equation*}
in $\Omega $, open neighborhood of the origin in $\mathbb{R}^{n}$, where
$r_{i}$, $i=2,\dots, n$, $k_{i-m+1}$ and $\ell_{i-m}$, $i=m+3,\dots, n$, are positive
integers such that $r_{2}< \cdots<r_{n}$, $k_{3}< \cdots<k_{n-m}$, $\ell_{3}< \cdots<\ell_{n-m}$ and
$r_{i}k_{i-m} < \ell_{i-m}$ for every $i$, $i = m+3,\dots, n$, then
the point $(0,e_{n})\in \Char(\widetilde{P}_{m,n})$ does not belong to  $WF_{\frac{r_{n}+r_{n-m}k_{n-m}}{r_{2}}}(u)$.
\end{remark}
%

%

%
\end{document}